\documentclass[12pt]{article}
\usepackage{makeidx}
\usepackage{amssymb}
 \usepackage{latexsym}
\usepackage{tikz}
\usepackage{pgfplots}
\usepackage{pgfkeys}%
\usetikzlibrary{plotmarks}
\usetikzlibrary{calc}

\usepackage{epsfig}
\usepackage{array}
\usepackage[ansinew]{inputenc}
\usepackage{amsfonts,amsmath}
\usepackage{latexsym}
\usepackage{graphics,color,graphicx,shortvrb}

 \usepackage{tabularx}
\newtheorem{Theorem}{Theorem}[section]
\newtheorem{Definition}[Theorem]{Definition}
\newtheorem{Proposition}[Theorem]{Proposition}

\definecolor{ProcessBlue}{cmyk}{1,0,0,0.25}
\definecolor{Black}{cmyk}{0,0,0,1}
\definecolor{Red}{cmyk}{0,1,1,0}
\definecolor{Green}{cmyk}{0.9,0,1,0}
\definecolor{Orange}{cmyk}{0,0.61,0.87,0.1}
\definecolor{Fuchsia}{cmyk}{0.47,0.91,0,0.06}
\definecolor{PineGreen}{cmyk}{0.92,0,0.59,0.25}

\usepackage[plainpages=false,hyperindex=false]{hyperref}

\begin{document}
\sloppy
 
 \begin{center} 
{\Large \bf Retro-Prospective Differential Inclusions and their 
Control by the Differential Connection Tensors of their 
Evolutions: The trendometer}\\ \mbox{}\\ Jean-Pierre 
Aubin\footnote{VIMADES (Viabilité, Marchés, Automatique, 
Décisions), 14, rue Domat, 75005, Paris, France \\
 aubin.jp@gmail.com, \url{http://vimades.com/aubin/}},\footnote{\textbf{Acknowledgments} 
\emph{This work was partially supported by the Commission of the 
European Communities under the 7th Framework Programme Marie 
Curie Initial Training Network (FP7-PEOPLE-2010-ITN),  project 
SADCO, contract number 264735.}} 
\\ \today


\end{center}
 
\begin{center} 
\textbf{Abstract}
\end{center} 

\emph{This study is motivated by two different, yet, connected, 
motivations. The first one follows the observation that the 
classical definition of derivatives involves prospective (or 
forward) difference quotients, not known whenever the time is 
directed,   at least at the macroscopic level.  Actually,  the 
available and known derivatives are retrospective (or backward). 
They coïncide whenever the functions are differentiable in the 
classical sense, but not in the case of non smooth maps, 
single-valued or set-valued. The later ones are used in 
differential inclusions (and thus, in uncertain control systems) 
governing evolutions in function of time and state. We follow the 
plea of some physicists for taking also into account the 
retrospective derivatives to study prospective evolutions in 
function of time, state and retrospective derivatives, a 
particular, but specific, example of historical of ``path 
dependent'' evolutionary systems. This is even more crucial in 
life sciences, in the absence of experimentation of uncertain 
evolutionary systems. The second motivation emerged from the 
study of networks with junctions (cross-roads in traffic 
networks, synapses in neural networks, banks in financial 
networks, etc.), an important feature of ``complex systems''.  At 
each junction, the velocities of the incoming (retrospective) and 
outgoing (prospective) evolutions are confronted. One measure of 
this confrontation (``jerkiness'') is provided by the product of 
the retrospective and prospective velocities, negative in 
``inhibitory'' junctions, positive for ``excitatory'' ones, for 
instance. This leads to the introduction of the ``differential 
connection tensor'' of two evolutions, defined as the tensor 
product of retrospective and prospective derivatives, which can 
be used for controlling evolutionary systems governing the 
evolutions through networks with junctions.} 

\mbox{}

\textbf{Mathematics Subject Classification}: 34A60,90B10, 90B20, 
90B99, 93C10, 93C30,93C99,
 
\textbf{Keywords} Transport, networks, junction, impulse, 
viability, traffic control, jam, celerity, monad

\section*{Motivations}

There are two different motivations of this study.

\subsection*{Retrospective-Prospective Differential Inclusions}
The first motivation follows the plea of \glossary{Galperin (Efim 
A.) [1934-]} \emph{Efim Galperin} in 
\cite[Galperin]{Galperin1,Galperin2,Galperin3,Galperin4}  for 
using  ``retrospective'' derivatives\footnote{For evolutions 
(fonctions of one variable, retrospective derivatives are 
derivatives from the left and prospective derivatives are 
derivatives from the right. For fonctions of several variable, 
there is no longer left and right, so retrospective (or 
backward), prospective (or forward), are used instead. }  instead 
of  ``prospective'' derivatives, universally chosen since their 
introduction by \emph{Newton} and \emph{Leibniz}, at a time when 
physics became predictive and deterministic: the ``prospective 
derivatives''  $\overrightarrow{D}x(t)$ being (more or less weak) 
limits of \emph{prospective (future) difference quotients} (on 
positive durations $h>0$)    
$\displaystyle{\overrightarrow{\nabla }_{h}x(t)  := 
\frac{x(t+h)-x(t)}{h}}$    are ``physically non-existent'', 
because they \textbf{are not yet known} at time $t$. Whereas the 
\emph{the retrospective (past) difference quotients 
$\displaystyle{\overleftarrow{\nabla }_{h}x(t) := 
\frac{x(t)-x(t-h)}{h}}$    \textbf{may be known}} for some 
positive durations and should be taken into account\footnote{This 
has been pointed out by     \emph{Jiri Buquoy}, who in 1812, 
formulated the equation of motion of a body with variable mass, 
which retained only the attention of \glossary{Poisson (Siméon 
Denis) [1781-1840]} \emph{Poisson} before being almost 
forgotten.  See  \cite[Jiri Buquoy]{Buquoy}, 
\cite[Mestschersky]{Mestschersky} and 
\cite[Levi-Civita]{Levi-Civita} among the precursors in this 
area.}. 

This is an inescapable issue in life sciences,  since the 
evolutionary engines evolve with time, under contingent and/or 
tychastic uncertainties and, in most cases, cannot be re-created 
(at least, for the time, since synthetic biology deals with this 
issue\footnote{See for instance \cite[Danchin et 
al.]{DanchinAll11}.}).  Popper's recommandations are valid for 
physical sciences, where experimentation is possible and 
renewable.    However,  the quest of the \index{instant} 
\emph{instant} (temporal window with $0$ duration) has not yet 
been experimentally created (the smallest measured duration is of 
the order of the yoctosecond ($10^{-24}$)). Furthermore, our 
brains deal with observations which are not instantaneous, but, 
in the best case, are perceived after a positive transmittal 
duration. 

For overcoming this difficulty, Fermat, Newton, Leibniz and 
billions of human brains have invented \emph{instants} and passed 
to the limit when duration of temporal windows goes to $0$ to 
reach such an instant. This is actually an 
approximation\footnote{actually, an inductive approximation, 
whereas (deductive) application refers to approximate derivatives 
of the idealized world by difference quotients, which are closer 
to the actual perception of our brains and capabilities of the 
digital computers.} of reality by clever mathematical 
constructions of objets belonging to an ever evolving ``cultural 
world''. Derivatives are not \emph{perceived}, but were 
\emph{invented}, simplifying reality by passing to the limit in a 
mathematical paradise.  

Therefore, for differentiable functions in the classical sense, 
the limits of retrospective and prospective difference quotients 
may coïncide when we pass to the limit. But this is no longer the 
case when evolutions are no longer differentiable in the 
classical sense, but derivatives may still exist for ``weaker'' 
limits, such as limits in the sense of distributions or graphical 
limits in set-valued analysis (see Section~18.9, p. 769, of 
\mbox{Viability Theory.  New Directions}, \cite[Aubin, Bayen \& 
Saint-Pierre]{absp}). Even if we restrict our analysis to 
Lipschitz functions, the \glossary{Rademacher (Hans Adolph) 
[1892-1969]} \emph{Rademacher}'s Theorem states that Lipschitz 
maps from one finite dimensional vector space to another one  are 
only \emph{almost everywhere} differentiable. Although small, the 
set of elements where there are not differentiable is interesting 
because Lipschitz have always set-valued graphical derivatives. 
Hence we have to make a detour by recalling what are meant 
retrospective and prospective graphical derivatives of maps as 
well as  set-valued maps and non differentiable (single-valued) 
maps. 
 
Therefore, we devote the first part of this study to a certain 
class of viable evolutions governed by functional (or 
history-dependent) differential inclusions  
$$x'(t) \; \in \;  G(t,x(t), \overleftarrow{D}x(t))$$ where 
$\overleftarrow{D}x(t)$ is the \emph{retrospective derivative} 
(or derivative from the left since, at this stage, we consider 
evolutions defined on $\mathbb{R}^{}$).  
Retrospective-prospective differential inclusions $x'(t) \; \in 
\;  G(t,x(t), \overleftarrow{D}x(t))$ describe \emph{predictions 
of evolutions based on the state \emph{and} on the known 
retrospective velocity at each chronological time}.  As delayed 
differential equations or inclusions, they are particular cases 
of \emph{functional} (or \emph{historical}, 
\emph{path-dependent}, etc.) differential equations\footnote{See 
\emph{Introduction to Functional Differential Equations}, 
\cite[Hale]{Hale},  \cite[Haddad]{hg81,hg81b,hg81c}, summarized 
in Chapter~12 of \emph{Viability Theory}, \cite[Aubin]{avt},  
\cite[Aubin \& Haddad]{aubhadclio,ah01hyb}, etc.}.  As for 
second-order differential equations, initial conditions 
$x(t_{0})$ at time $t_{0}$ must be provided, as well as 
(retrospective) initial velocities for selecting evolutions 
governed by retrospective-prospective differential equations.

\subsection*{Differential Connection Tensors in Networks}

The second motivation emerged from the study of propagation 
through ``junctions of a network'', such as cross-roads in road 
networks, banks in financial networks, synapses in neural 
network, etc. (see for instance \cite[Aubin]{TransReg-2013}).   

\subsubsection*{Neural Network : the Hebbian Rule} 

If we accept that in formal neuron networks,  ``(evolving) 
knowledge'' is coded as ``synaptic weights'' at each synapse, 
their collection defines a ``synaptic matrix'' which evolves, 
and, thus, becomes the ``state of the network''.  \emph{Donald 
Hebb} introduced in 1949   \glossary{Hebb (Donald) [1904-1985]}  
in   \emph{The Organization of Behavior}, \cite[Hebb]{Hebb}, the 
\emph{Hebbian learning rule} prescribing that  the velocity of 
the synaptic matric is proportional to the \index{tensor product} 
\emph{tensor product}\footnote{Recall that the tensor product $p 
\otimes q$ of two vectors $p  := (p_{i} )_{i} \in 
\mathbb{R}^{\ell} $ and $q  :=(q_{j} )_{j} \in \mathbb{R}^{\ell} 
$ is the rank one linear operator $$p  \otimes q \in 
\mathcal{L}(\mathbb{R}^{\ell},\mathbb{R}^{\ell}) : x \mapsto 
\left\langle p,x \right\rangle q $$ the entries of which (in the 
canonical basis) are equal to $(p_{i} q_{j}) _{i,j}$.} of the 
``presynaptic activity'' and ``postsynaptic activity'' described 
by the propagation of nervous influx in the neurons.

Hence, denoting the synaptic matrix $W$ of synaptic weights, the 
basic question was to minimize a ``matrix function'' $W \in 
\mathcal{L}(X,X)\mapsto E(Wx)$ where $x \in X:= 
\mathbb{R}^{\ell}$ and $E: X \mapsto \mathbb{R}^{}$ a 
differentiable function are given. Remembering\footnote{See 
Proposition~2.4.1, p. 37 and Chapter~2 of  \mbox{Neural Networks 
and Qualitative Physics: a Viability Approach}, 
\cite[Aubin]{a92ia}.} that the gradient with respect to $W$ is 
equal to the tensor product $E'(Wx) \otimes x$, the gradient 
method leads to a differential equation of the form 
\begin{equation} \label{e:}   
W'(t) \; = \; - \alpha E'(W(t)x) \otimes x 
\end{equation}
which governs the evolution of the synaptic matrix (the ``synapse 
$x$ is fixed and does not evolve). 

\subsubsection*{Differential Connection Tensors}

However, we take into account  the evolution $t \mapsto x(t) \in 
X $ of the propagation in networks (such as the propagation 
nervous influx,  traffic,  financial product, etc.). If the 
evolution is Lipschitz,  retrospective and prospective 
derivatives exist at all times, so that we can define the tensor 
product $\overleftarrow{D}x(t)\otimes \overrightarrow{D}x(t)$ of 
their retrospective and prospective velocities: we shall call it 
the \index{} \emph{differential connection tensor} of the 
evolution $x(\cdot)$ at time $t$.

It \emph{plays the role of a \index{trendometer}  
\emph{``trendometer''} measuring the \index{monotonicity 
reversal} \index{trend reversal} \emph{trend reversal} (or 
\emph{monotonicity reversals}) at junctions}:  the differential 
connection tensor describes the \emph{trend reversal} between the 
retrospective and prospective   trends  when they are strictly 
negative, the \emph{monotonicity congruence} when they are 
strictly positive and the \emph{inactivity} they  vanish. In 
neural networks, for instance, \emph{this an inhibitory effect or 
trend reversal in the first case, an excitatory or trend 
congruence in the second case, and inactivity of a synapse: one 
at least of the propagation of the nervous influx stops.}   The 
absolute value of this product measures in some sense the 
\emph{jerkiness of the trend reversal} at a junction of the 
network.

We are thus tempted to control (pilot, regulate, etc.) the 
evolution of  propagation in the network   governed by a system 
\begin{equation} \label{e:synpropagactivity}   
x' (t) \; = \; g (x (t), u (t)) \;\mbox{\rm where}\; u (t) \; \in 
\;  U (\overleftarrow{D}x(t)\otimes \overrightarrow{D}x(t))
\end{equation}
controlled by differential connection tensors at junctions of the 
network. We recall that the evolutions governed by (Marchaud) 
controlled systems are Lipschitz under the standard assumption, 
\emph{but not necessarily differentiable}. For example, in order 
to  govern the viability of the propagation in terms of the 
inhibitory, excitatory and stopping behavior at  the junctions of 
the network,  some constraints are imposed on the evolution of 
the differential connection tensors.  Examples of 
retrospective-retrospective differential equations are provided 
by tracking or controlling differential connection tensors of the 
evolutions requiring that evolutions governed by differential 
equations $x'(t)=f(t,x(t))$ satisfy contraints of the form 
$\overleftarrow{D}x(t) \otimes \overrightarrow{D}x(t) \in 
C(t,x(t))$. These control  systems are examples of 
retrospective-prospective differential inclusions.

These considerations extend to  ``multiple synapses'' when we 
associate with each subset $S$  of branches $j$ meeting at a 
junction  the tensor products $\otimes_{j \in S}x_{j}'(t) $ of 
the velocities at the junction\footnote{See \cite[Aubin  \& 
Burnod]{ab99nng} and the literature on $\Sigma-\Pi$ neural 
systems, Section~12.2 of \mbox{Viability Theory.  New 
Directions}, \cite[Aubin, Bayen \& Saint-Pierre]{absp}, 
\emph{Analyse qualitative}, \cite[Dordan]{do93liv}, as well as 
\cite[Aubin]{a92ia, aconcom98,a01ctfcr,AUB-Leit09,TransReg-2013}, 
\cite[Vinogradova]{vino} and the literature on the regulation of 
networks.}.

 \subsection*{Organization of the Study}

Section~\ref{s:Remarks}, p. \pageref{s:Remarks}, 
\emph{Retrospective-Prospective Differential Inclusions},  
defines retrospective and prospective (graphical) derivatives of 
tubes and evolutions, their \emph{differential connection tensor} 
(Definition~\ref{d:DifConnTensorTube}, 
p.\pageref{d:DifConnTensorTube}). They are the ingredients   for 
introducing retrospective-prospective differential inclusions. 
The Viability Theorem (Theorem~\ref{t:RPdifInc}, 
p.\pageref{t:RPdifInc}) is  adapted for characterizing viable 
tubes under such differential inclusions using characterizations 
linking the retrospective and prospective derivatives of the 
tube. When these conditions are not satisfied, we restore the 
viability by introducing \emph{retrospective-prospective 
viability kernel} of the tube under the retrospective-prospective 
differential inclusion (Subsection~\ref{s:PRKernel}, p. 
\pageref{s:PRKernel}).

Section~\ref{s:ControlDCT}, p. \pageref{s:ControlDCT}, 
\emph{Control by Differential Connection Tensors}, studies the 
regulation of viable evolutions on tubes by imposing constraints 
on their differential connection tensors.
 
Section~\ref{s:Illustrations}, p. \pageref{s:Illustrations}, 
\emph{Illustrations}, provides examples of differential 
connection tensors of vector evolutions in the framework of 
``technical analysis'' of the forty prices series of the CAC 40 
stock market index\footnote{See Chapter~2 of \emph{Tychastic 
Viability Measure of Risk Eradication. A Viabilist Portfolio 
Performance and Insurance Approach}, \cite[Aubin, Chen Luxi \& 
Dordan]{ACD-ALIM}, for a more detailed study.}. 

Section~\ref{s:Remarks}, p. \pageref{s:Remarks},  \emph{Other 
Examples of Differential Connection Tensors}, defines 
differential connection tensors of set-valued maps 
(Subsection~\ref{s:PRDerivatives}, p. \pageref{s:PRDerivatives}, 
\emph{Prospective and Retrospective Derivatives of Set-Valued 
Maps}) and gathers some other classes  differential connection 
tensors than the ones of the evolutions $t \mapsto x(t)$ or tubes 
$t  \leadsto K(t)$ from $\mathbb{R}^{}$ to $\mathbb{R}^{\ell}$, 
which provided the first source of motivations for studying 
differential connection tensors. Other specific examples are the 
\emph{differential connection tensors of numerical functions} $V 
: \mathbb{R}^{\ell} \mapsto \mathbb{R}^{}$ 
(Subsection~\ref{s:PREpider}, p. \pageref{s:PREpider}), and 
\emph{tangential connection tensors} of retrospective and 
prospective tangents (Subsection~\ref{s:PRCones}, p. 
\pageref{s:PRCones}). These issues are the topics of forthcoming 
studies.

\section{Retrospective-Prospective Differential Inclusions} \label{s:PRTubes}

\subsection{Prospective and Retrospective Derivatives of Tubes and 
Evolutions} 

A tube is the nickname of a set-valued map $K: t \in 
\mathbb{R}^{}  \leadsto K(t) \subset X$. Since there are 
only\footnote{Actually, there is a third one, $0$, where 
$\overleftarrow{D}K(t,x)(0)$ and $\overrightarrow{D}K(t,x)(0)$ 
are the retrospective and prospective tangent cones studied in 
Section~\ref{s:PRCones}, p. \pageref{s:PRCones}.} two directions 
$+1$ and $-1$ in $\mathbb{R}^{}$,  the prospective (left) and 
retrospective (right) derivatives of a tube $K$ at a point 
$(t,x)$ of its graph are defined by 
\begin{equation} \label{e:} \left\{ \begin{array}{ll}
\displaystyle{v \in \overrightarrow{D}K(t,x) \;\mbox{\rm if and 
only if }\; \liminf_{h \rightarrow 0+} d\left(v, 
\frac{K(t+h)-x}{h}\right) \; = \; 0 }\\
\displaystyle{v \in \overleftarrow{D}K(t,x) \;\mbox{\rm if and 
only if }\; \liminf_{h \rightarrow 0+} d\left(v, 
\frac{x-K(t-h)}{h}\right) \; = \; 0 }\\
\end{array} \right.  \end{equation}
(see Definition~\ref{d:PRderiv}, p.\pageref{d:PRderiv}. in the 
general case). 

\begin{Definition} 
\symbol{91}\textbf{Differential Connection Tensor of a 
Tube}\symbol{93} \label{d:DifConnTensorTube}\index{} The 
\index{differential connection tensor of a tube} 
\emph{differential connection tensor of a tube} $K(\cdot)$ at $x 
\in K(t) $ is defined by
\begin{equation} \label{e:ReverIndex}  
\forall \; \overleftarrow{v} \in \overleftarrow{D}K(t,x), \; 
\forall \; \overrightarrow{v} \in  \overrightarrow{D}K(t,x), \; 
\; \mathbf{a}_{K}(t,x)[\overleftarrow{v},\overrightarrow{v}] \; 
:= \;  \overleftarrow{v}\otimes\overrightarrow{v} 
 \end{equation}
\end{Definition}

In particular, an evolution $x(\cdot)$ is a single-valued tube 
defined by $K(t):=\{x(t)\}$, so that we can define their 
graphical prospective derivative $\overrightarrow{D}x(t)$ (from 
the right) and    retrospective derivatives 
$\overleftarrow{D}x(t)$ (from the left) respectively (see 
illustrations in Section~\ref{s:Illustrations}, p. 
\pageref{s:Illustrations}, \emph{Illustrations}).

\subsection{Retrospective-Prospective Differential Inclusions}

Recall that whenever an evolution $t \mapsto x(t)$ is viable on a 
neighborhood of $t_{0}$ on a tube $K(t)$, then 
$\overleftarrow{D}x(t_{0}) \; \in \;  
\overleftarrow{D}K(t_{0},t_{0}) $ and $\overrightarrow{D}x(t_{0}) 
\; \in \;  \overrightarrow{D}K(t_{0},t_{0}) $.

Since we know only retrospective derivatives,  forecasting   
future evolution can be governed  by prospective differential 
inclusion  $\overrightarrow{D}x(t) \in F(t,x(t))$ depending only 
on time and state, but also by the particular case of 
history-dependent evolutions $\overrightarrow{D}x(t) \in 
G(t,x(t), \overleftarrow{D}x(t))$ depending on  time, state and 
the retrospective derivatives. This could be the case for system 
controlling the differential connection tensors  of the 
evolutions, for instance (see Section~\ref{s:ControlDCT}, p. 
\pageref{s:ControlDCT}).

\begin{Theorem} 
\symbol{91}\textbf{Viability Theorem for 
Retrospective-Prospective  Differential Inclusions}\symbol{93} 
\label{t:RPdifInc}\index{} Let us assume that the map $(t,x,v) 
\in \mathbb{R}^{} \times X \times X  \leadsto G(t,x,v) \subset X$ 
is Marchaud (closed graph, convex valued and linear growth) and 
that the tube $t  \leadsto K(t)$ is closed. Then the ``tangential 
condition'' 

\begin{equation} \label{e:RPdifInc}   
\forall \; \overleftarrow{v} \in \overleftarrow{D}K(t,x), \; \; 
G(t,x,\overleftarrow{v}) \cap \overrightarrow{D}K(t,x) \; \ne  \; 
\emptyset
\end{equation} 
is equivalent to the ``viability property'': from any initial 
state $x_{0} \in K(t_{0})$ and initial retrospective velocity 
$\overleftarrow{v}_{0} \in \overleftarrow{D}K(t_{0},x_{0})$, 
there exists at least one evolution $x(\cdot)$ governed by the 
retrospective-prospective differential inclusion 
$\overrightarrow{D}x(t) \in G(t,x(t),\overleftarrow{D}x(t))   $ 
satisfying  $x(t_{0})=x_{0}$ and $\overleftarrow{D}x(t_{0})= 
\overleftarrow{v}_{0}$ and viable in the tube $K(\cdot)$.  
\end{Theorem} 
 
\textbf{Proof} --- \hspace{ 2 mm} The proof is an adaptation of 
the proof of the viability Theorem~19.4.2, p. 782, based on 
Theorems~11.2.7, p. 447, and 19.3.3, p. 777, of  \emph{Viability 
Theory.  New Directions}, \cite[Aubin, Bayen \&  
Saint-Pierre]{absp}. We just indicate the modifications to be 
made.

We  construct approximate solutions   by modifying Euler's method 
to take into account the viability constraints, then deduce from 
available estimates that a subsequence of these solutions 
converges in some sense to a limit, and finally, check that this 
limit is a viable solution to the retrospective-prospective  
differential inclusion ($\overrightarrow{D}x(t) \in 
G(t,x(t),\overleftarrow{D}x(t))  $.

\begin{enumerate}

\item By assumption, there exists $r>0$ such that the 
neighborhood  $\mathcal{K}_{r} := \mbox{\rm Graph}(K)  \cap 
(t_{0},x_{0})+r([-1,+1] )\times B$ of the initial condition 
$(t_{0},x_{0})$  is compact. Since $G$ is Marchaud, the set
\begin{displaymath}
\mathcal{C}_{r}  \;:=  \; \left\{ F(t,x, \overleftarrow{v}) 
\right\} + B, \;  \; \mbox{\rm and} \; T  \;:= \; 
r/\|\mathcal{C}_{r}\|
\end{displaymath}
is also compact. We next associate with any $h$ the Euler 
approximation

\begin{equation} 
v_{j}^{h} \; := \; \frac{x_{j+1}^{h}-x_{j}^{h}}{h}  \; \in \; 
G(jh,x_{j}^{h}, v_{j-1}^{h}) \;\mbox{\rm where}\;  v_{j-1}^{h} \; 
:= \;  \frac{x_{j}^{h}-x_{j-1}^{h}}{h}  
\end{equation}
starting from $(t_{0}, x_{0}, \overleftarrow{v}_{0})$.

\item Theorems~11.2.7, p. 447 of \cite[Aubin, Bayen \&  
Saint-Pierre]{absp} implies that for all $\varepsilon >0$, 
\begin{equation} \label{eq-apprexplschviab2} \left\{ 
\begin{array}{l} \exists \; \eta (\varepsilon) >0 \; \mbox{such 
that} \; \forall \; (t,x) \in \mathcal{K}_{r}, \; \forall h \in 
\left[0,\eta (\varepsilon) \right], \\ x_{j}^{h}+ 
hG(jh,x_{j}^{h}, v_{j-1}^{h}) \; \in \; K(jh,,x_{j}^{h}) 
+\varepsilon B
\end{array} \right.
\end{equation}
Since \begin{displaymath} \|x_{j}^{h}-x_{0}\| \;\leq  \; 
\sum_{i=0}^{i=j-1}\|x_{i+1}^{h}-x_{i}^{h}\| \; \leq  \; 
\sum_{i=0}^{i=J^{h}-1}h \left\| v_{j}^{h}\right\| \; \leq \; 
\|\mathcal{C}_{k}\|
\end{displaymath}
the discrete evolution is viable in $\mathcal{K}_{r}$   on the 
interval $[0,T]$. Denoting by $x^{h}$, $\overleftarrow{v}^{h}$ 
and $\overrightarrow{v}_{h}$ the linear interpolations of the 
sequences $x^{h}_{j}$, $\overleftarrow{v}^{h}_{j}$ and 
$\overrightarrow{v}^{h}_{j}$, we infer that there exists a 
constant $\alpha>0$ such that

\begin{equation} \label{e:} \left\{ \begin{array}{l} 
\displaystyle{ (t^{h},x^{h}, 
\overleftarrow{v}^{h},\overrightarrow{v}) \; \in \; \mbox{\rm 
Graph}(G) + \varepsilon \alpha}\\
\displaystyle{(t^{h},x^{h}) \; \in \; \mbox{\rm Graph}(K) + \varepsilon \alpha}\\
\end{array} \right. 
\end{equation}
and that there exists a constant $\beta >0$ such that   the 
\textit{a priori} estimates

\begin{equation} \label{xA445} 
\max(\|x^{h}\|_{\infty}, \|\overleftarrow{\nabla} ^{h}  
x^{h}\|_{\infty},  |\overrightarrow{\nabla} ^{h}  
x^{h}\|_{\infty})  \; \leq \;  \beta
\end{equation}
are satisfied.
   
\item   They imply the \emph{a priori} estimates   of the 
Convergence Theorem~19.3.3, p. 777, of    \cite[Aubin, Bayen \&  
Saint-Pierre]{absp}, which states   the limit of a converging 
subsequence is a solution to the retrospective-prospective 
differential inclusion, viable in $\mbox{\rm Graph}(K)$. \hfill 
$\;\; \blacksquare$ \vspace{ 5 mm}
\end{enumerate}

 
\subsection{Retrospective-Prospective Viability Kernels} 
\label{s:PRKernel} 


Naturally, the ``tangential assumption'' (\ref{e:RPdifInc}), p. 
\pageref{e:RPdifInc}, is not necessarily satisfied so that we 
have to adapt the concept of viability kernel to the 
retrospective-prospective case.

\begin{Definition} 
\symbol{91}\textbf{Retrospective-Prospective Viability Kernel of 
a Tube}\symbol{93}\label{}\index{} The viability kernel of the 
tube $K(\cdot)$ is the set of initial conditions $(t_{0},x_{0}, 
\overleftarrow{v}_{0}) \in \mathbb{R}^{} \times K(t_{0}) \times 
\overleftarrow{D}K(t_{0},x_{0})$ from which starts at least one 
viable evolution $t \mapsto x(t) \in K(t)$ to the 
retrospective-prospective differential inclusion in the sense that

\begin{equation} \label{e:} \left\{ \begin{array}{ll}
 (i) & \overrightarrow{D}x(t) \; \in \;  G(t,x(t),\overleftarrow{D}x(t))\\ 
 (ii) & \overleftarrow{D}x(t) \; \in \; \overleftarrow{D}K(t,x(t)) 
\;\mbox{\rm and}\; \overrightarrow{D}x(t) \; \in 
\overrightarrow{D}K(t,x(t))\; 
\\
\end{array} \right. 
\end{equation}
 
\end{Definition}

We provide a viability characterization of 
retrospective-prospective viability kernel tubes:

\begin{Proposition} 
\symbol{91}\textbf{Viability Characterization of 
Retrospective-Prospective Viability 
Kernel}\symbol{93}\label{}\index{} Let us consider the control 
system

\begin{equation} \label{e:RPcontSysy} \left\{ \begin{array}{ll} 
(i) & \tau'(t) \; = \; 1\\ 
(ii) &  x'(t) \; \in \; G(\tau(t),x(t), \overleftarrow{v}(t))\\
(iii) & \|\overleftarrow{v}'(t)\|  \; \leq \; c \; \|G(t,x,\overleftarrow{v})\|\\
 & \mbox{\rm where}\; \overleftarrow{v}(t) \; \in \;  
\overline{\mbox{\rm co}} (\overleftarrow{D}K(\tau(t),x(t)))
\end{array} \right. 
\end{equation}
Then the viability kernel of the graph $\mbox{\rm 
Graph}(DK(\cdot))$ of the derivative tube $K(\cdot)$ coincides 
with the retrospective-prospective viability kernel of the tube. 
\end{Proposition}

\textbf{Proof} --- \hspace{ 2 mm} The viability kernel of the 
control system (\ref{e:RPcontSysy}), p. \pageref{e:RPcontSysy} is 
the set of initial triple $(t_{0}, x_{0}, \overleftarrow{v}_{0})$ 
such that $x_{0} \in K(t_{0})$ and $\overleftarrow{v}_{0} \in 
\overleftarrow{D}K(t_{0},x_{0})$ from which starts an evolution  
$t \mapsto (t_{0}+t, x(t),\overleftarrow{v}(t))$ of the control 
system such that  $x(t) \in K(\tau(t))$ and $\overleftarrow{v}(t) 
\; \in \; \overline{\mbox{\rm co}} 
(\overleftarrow{D}K(\tau(t),x(t)))$. Setting 
$x_{\star}(t):=x(t-t_{0})$ and  
$\overleftarrow{v}_{\star}(t):=\overleftarrow{v}(t-t_{0})$, we 
observe that $x_{\star}(t) \in 
G(t,x_{\star}(t),\overleftarrow{v}_{\star}(t))$, 
$\overleftarrow{v}_{\star}(t) \in 
\overleftarrow{D}K(t,x_{\star}(t))$ and $x_{\star}(t) \in K(t)$. 
We thus infer that $\overrightarrow{D}x_{\star}(t) \in 
\overrightarrow{D}K(t,x_{\star}(t))$. Since $x(t)$ is viable in 
the tube, we also infer that $\overleftarrow{D}x(t)$ actually 
belongs to $\overleftarrow{D}K(t,x(t))$. Hence $(t_{0}, x_{0}, 
\overleftarrow{v}_{0})$ belongs to the retrospective-prospective 
viability kernel of the tube $K(\cdot)$. \hfill  $\;\; 
\blacksquare$ \vspace{ 5 mm}

Therefore, it remains to provide sufficient conditions for the 
viability kernel of the graph of $K(\cdot)$ under the control 
system is Marchaud. 

\begin{Theorem} 
\symbol{91}\textbf{Properties of the Retrospective-Prospective 
Viability Kernel}\symbol{93}\label{} Let us assume that the 
set-valued map $G: (t,x, \overleftarrow{v})  \leadsto G(t,x, 
\overleftarrow{v})$ is Marchaud. Then the 
retrospective-prospective viability kernel of the tube $K(\cdot)$ 
under the $\overleftarrow{D}x(t) \in G(t,x(t), 
\overleftarrow{D}x(t))$ is closed and inherits all properties of 
viability kernels.
 \end{Theorem}

\section{Control by Differential Connection Tensors}
\label{s:ControlDCT} 


We study the tracking  at each date $t$ of the differential 
connection tensor $\overleftarrow{D}x(t)\otimes 
\overrightarrow{D}x(t)$ of evolutions  governed by a differential 
inclusion $x'(t) \in F(t,x(t)) $.  

For that purpose, we introduce a connection map $(t,x)  \leadsto 
C(t,x)  \;\subset\; \mathcal{L}(X,X)$. We are looking for 
evolutions $x(\cdot)$ governed by the differential inclusion 
satisfying the constraints on the differential connection tensors

\begin{equation} \label{e:}   
\forall \; t \geq 0, \; \; \overleftarrow{D}x(t)\otimes 
\overrightarrow{D}x(t) \; \in \; C(t,x(t))
\end{equation}

This is a problem analogous to the search of the slow evolutions 
governed by control systems (solutions governed by controls of 
the regulation map with minimal norm): see \cite[Aubin  \& 
Frankowska]{af85heav} or Theorem~6.6.3, p. 229, of 
\cite[\emph{Viability Theory}]{avt}.

We follow the same strategy 
by introducing   the set-valued map $ G$ defined by

\begin{equation} \label{e:TRregRev}   
G(t,x, \overleftarrow{v}) \; := \; \left\{w \in F(t,x) \; \mbox{ 
such that} \; \overleftarrow{v} \otimes w \; \in \; C(t,x) 
  \right\}
\end{equation}

\begin{Theorem} 
\symbol{91}\textbf{Control of Differential Connection 
Tensors}\symbol{93} \label{}\index{} We assume that $F$ is 
Marchaud, that the tube $K(\cdot)$ is closed and that

\begin{equation} \label{e:} \left\{ \begin{array}{ll} 
(i)&  \mbox{\rm the graph of $(t,x)  \leadsto C(t,x) \subset 
\mathcal{L}(X,X)$ is closed and its images 
are convex}\\
(ii) & \forall \; (t,x) \in \mbox{\rm Graph}(K), \; \; \forall \; 
\overleftarrow{v} \in \overleftarrow{D}K(t,x), \; \; \exists \;  
 w \in F(t,x) \in \overrightarrow{D}K(t,x)  \\ &  \mbox{ such that} 
\; \; \overleftarrow{v} \otimes w \; \in \; C(t,x)
\end{array} \right. 
\end{equation}
 
For any $t_{0}$, for any $x_{0} \in K(t_{0})$, for any 
$\overleftarrow{v}_{0} \in \overleftarrow{D}K(t_{0},x_{0})$, 
there exists at least an evolution $x(\cdot)$ governed by the 
differential inclusion $x'(t) \in F(t,x(t))$ starting at $x_{0}$ 
viable in the tube $K(\cdot)$ such that $\overleftarrow{v}_{0} 
\otimes \overrightarrow{D}x(t_{0}) \in C(t_{0},x_{0})$ and 
satisfying the differential connection tensor constraints
\begin{equation} \label{e:}   
\forall \; t \geq t_{0}, \; \; \overleftarrow{D}x(t) \otimes 
\overrightarrow{D}x(t) \; \in \;  C(t,x(t))
\end{equation}
and the retrospective-prospective viability property
\begin{equation} \label{e:}   
\forall \; t \geq t_{0}, \; \; \overleftarrow{D}x(t) \otimes 
\overrightarrow{D}x(t) \; \in \;  \overleftarrow{D}K(t,x(t)) 
\otimes  \overrightarrow{D}K(t,x(t)) 
\end{equation}

\end{Theorem}

{\bf Proof} --- \hspace{ 2 mm} The set-valued map $G$   satisfies 
the assumptions of Theorem~\ref{t:RPdifInc}, 
p.\pageref{t:RPdifInc}, in such a way that there exists one 
evolution $x(\cdot)$ governed by $\overrightarrow{D}x(t) \in 
G(t,x(t), \overleftarrow{D}x(t))$ viable in the tube $K(\cdot)$.  
Therefore, $\overrightarrow{D}x(t) \; \in \; 
\overrightarrow{D}K(t,x(t))$ for all $t \geq t_{0}$. 
Consequently, 

\begin{equation} \label{e:}   
\overleftarrow{D}x(t) \otimes \overrightarrow{D}x(t) \in C(t,x(t))
\end{equation}
and since the evolution is viable in the tube $K(\cdot)$, that

\begin{displaymath}  
\overleftarrow{D}x(t) \; \in \;  \overleftarrow{D}K(t,x(t)) 
\;\mbox{\rm and}\; \overrightarrow{D}x(t) \; \in \; 
\overrightarrow{D}K(t,x(t))
\end{displaymath}
The theorem ensues.   \hfill $\;\; \blacksquare$ \vspace{ 5 mm}
  
For instance, we can choose

\begin{equation} \label{e:}   
C(t,x,\overleftarrow{v}) \; := \; \left\{\overrightarrow{v} \; 
\mbox{ such that} \; \sup_{w \in F(t,x)} \sup_{(i,j)} 
\overleftarrow{v}_{i} (\overrightarrow{v}_{j} -w_{j})  \; \leq \; 
0   \right\}
\end{equation}

In other words, the  entries 
$\overleftarrow{v}_{i}\overrightarrow{v}_{j} $ minimize the 
entries $\overleftarrow{v}_{i} w_{j}$ of the differential 
connection tensors when the velocities $w \in F(t,x) $.  

Proposition~6.5.4, p. 226, of \emph{Set-valued analysis}, 
\cite[Aubin \& Frankowska]{af90sva}, implies that the connection 
constraint map has a closed  graph and  convex values whenever 
the set-valued map $F$ is lower semicontinuous with convex 
compact images.
 
We could as well requires that the entries of the differential 
connection tensor maximize the  entries 
$\overleftarrow{v}_{i}\overrightarrow{v}_{j} $ minimize the 
entries $\overleftarrow{v}_{i} w_{j}$ of the differential 
connection tensors when the velocities $w \in F(t,x) $ or that 
for some pairs $(i,j)$, the entries 
$\overleftarrow{v}_{i}\overrightarrow{v}_{j} $ minimize 
$\overleftarrow{v}_{i}w_{j} $ and for the other pairs, that they 
maximize $\overleftarrow{v}_{i}w_{j} $ when the velocities $w \in 
F(t,x) $.
 
 
\section{Illustrations} \label{s:Illustrations}
 
The question arises whether it is possible to detect the 
connection dates \emph{when the monotonicity of a series  of a 
family of temporal series is followed by the reverse (opposite) 
monotonicity of other series}, in order to detect the influence 
of each series on the dynamic behavior of other ones. When the 
two functions are the same, we obtain  their reversal dates when 
the series achieve their extrema. The \emph{differential 
connection tensor} measures the \index{jerkiness} 
\emph{jerkiness} between two functions, smooth or not smooth 
(temporal series) providing the trend reversal dates of the 
differential connection tensor.

This matrix plays for time series  a dynamic rôle  analogous to 
the static rôle played by the correlation matrix of a family of 
random variable measuring the covariance   entries between two 
random coefficients. In other words, we add in our analysis the 
dependence on \emph{random} events of  variables   their 
dependence on \emph{time}. 

The differential connection tensor softwares provides at each 
date  the coefficients of the differential connection tensor. 

We use the tensor trendometer for detecting the dynamic 
correlations between the forty price series of the  CAC 40. For 
instance, on August 6, 2010, the prices are displayed in the 
following figure

\begin{center} \includegraphics[width=.6\linewidth]{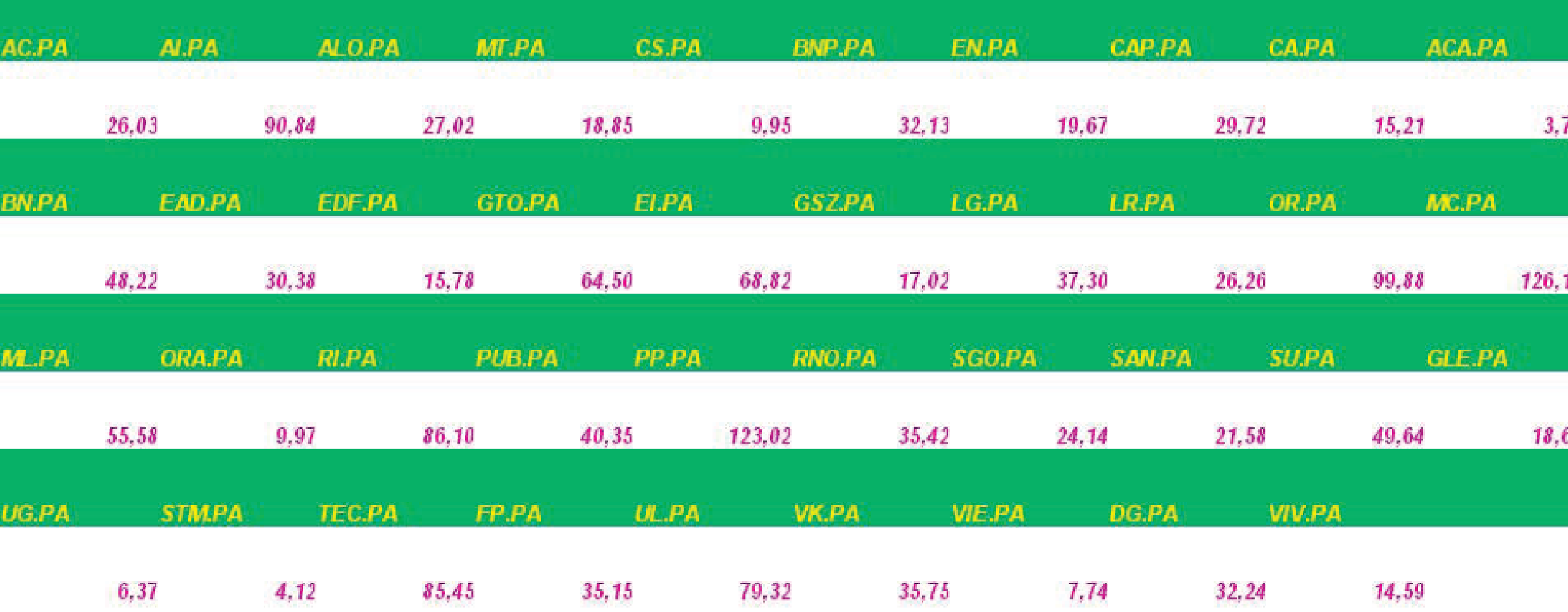}   
\end{center} 
 
At each date, it provides the $40 \times 40$ matrix displaying 
the qualitative jerkiness for each pair  of series when the trend 
of the first one is followed by the opposite trend of the second 
one. At each entry, the existence of a trend reversal by a 
circles:
  
\begin{center} \includegraphics[width=.6\linewidth]{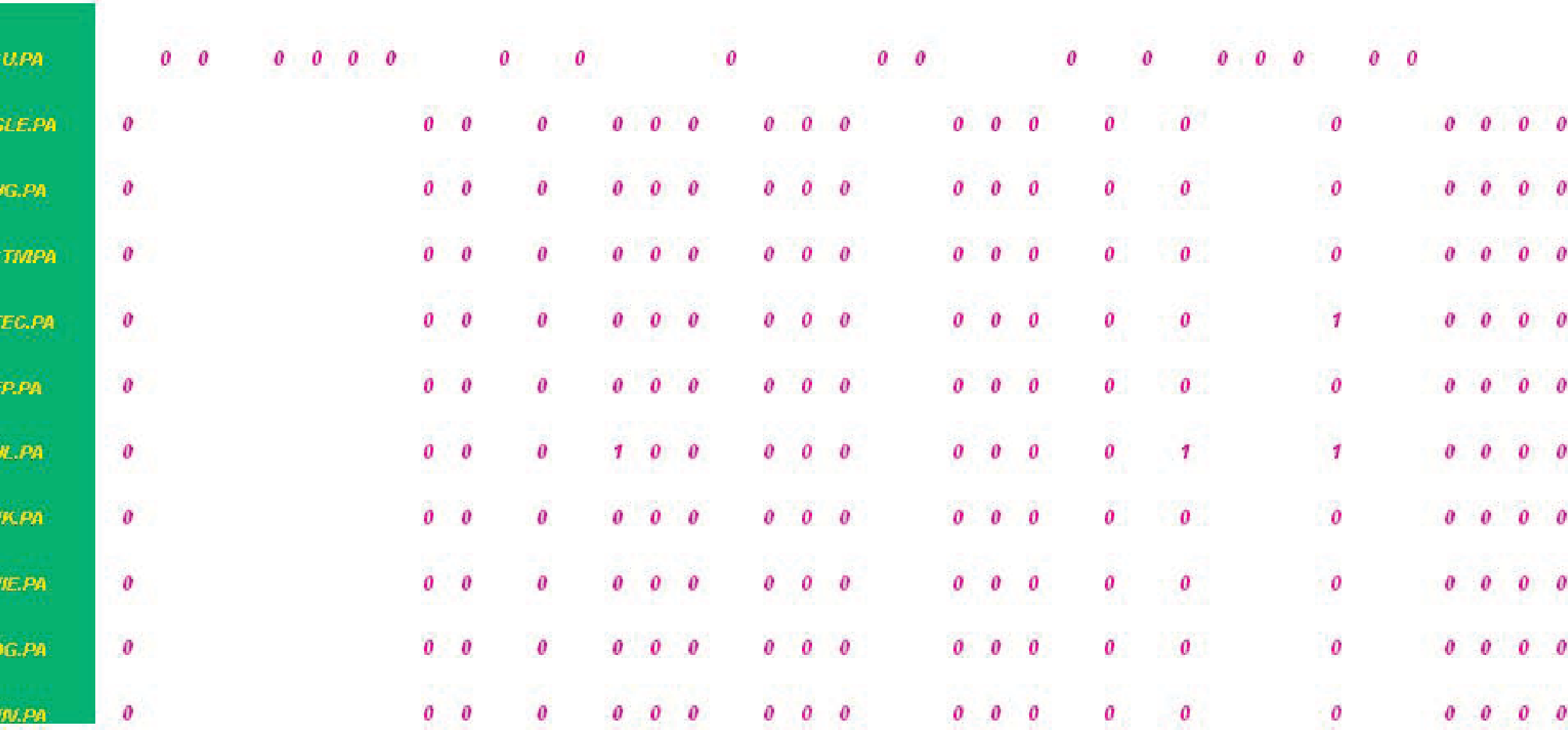}   
\end{center}

The quantitative version replaces the circles by the values of 
the jerkiness:

\begin{center} \includegraphics[width=.6\linewidth]{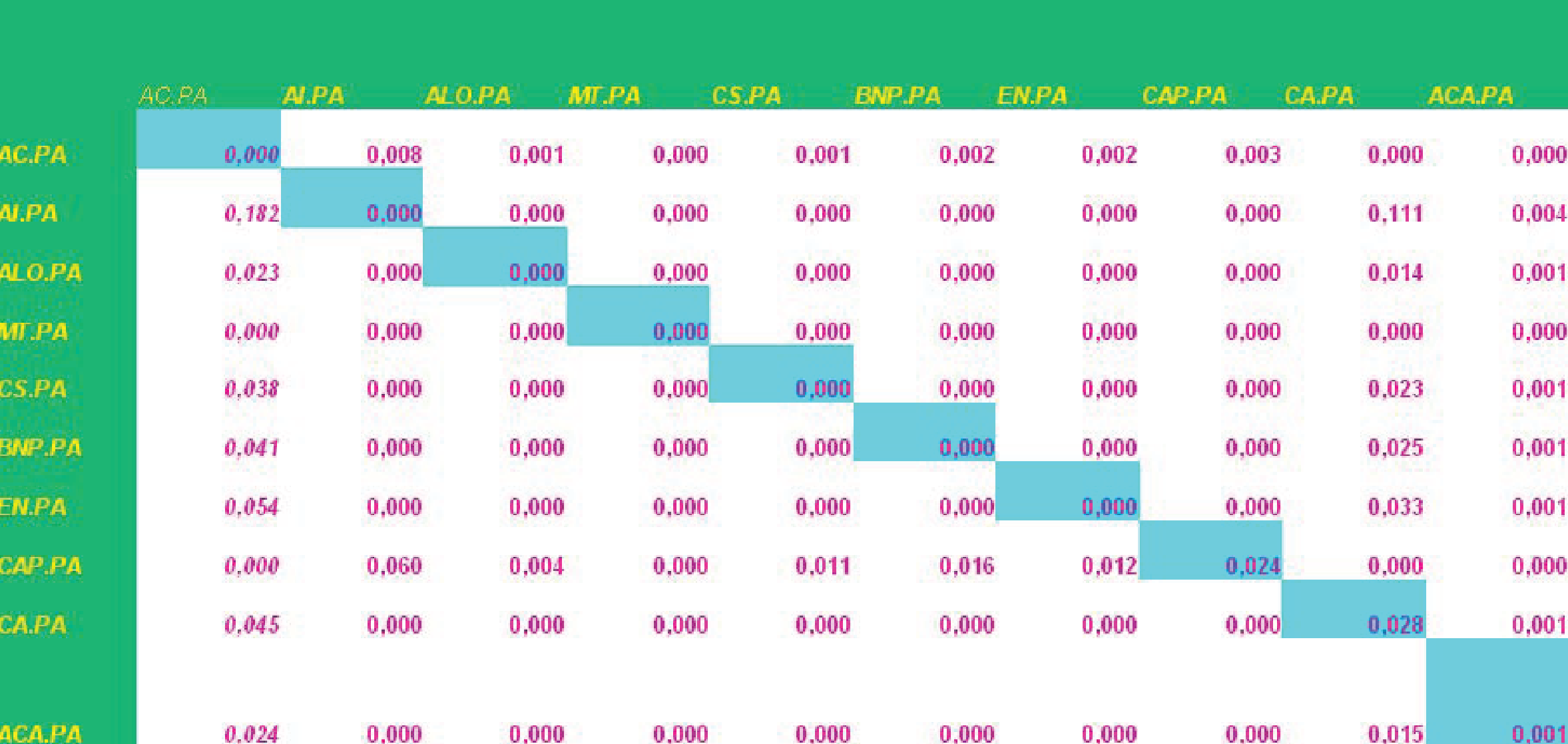}   
\end{center} 
 
In order to analyse further the evolutionary behavior of the CAC 
40, we present the analysis of the CAC 40 index only, but  over 
the period   from  du 03/01, 1990 to   09/25, 2013. The first 
figure displays the series of the CAC 40 indexes (closing 
prices). The vertical bars indicate the reversal dates and their 
height displays their jerkiness.

\begin{center} 
 \includegraphics[width=.6\linewidth]{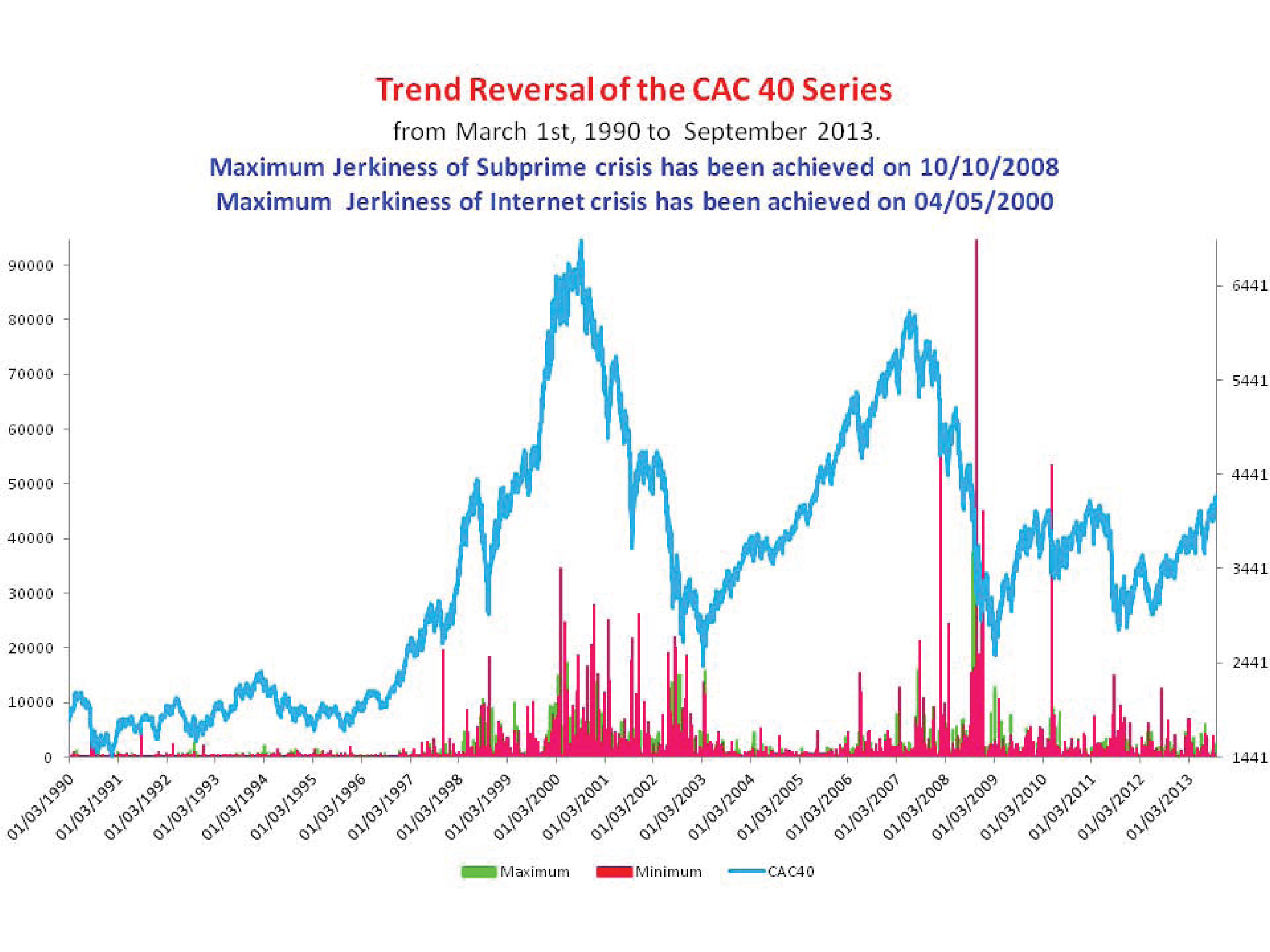} 
\end{center} 

The 2000 Internet crisis (around May 4, 2000) and the 2008 
``subprime'' crisis (around October 10, 2008)  are detected and 
measured by the trendometer:

\begin{center}\begin{tabular}{m{0.474\textwidth}m{0.474\textwidth}} 
\includegraphics[width=\linewidth]{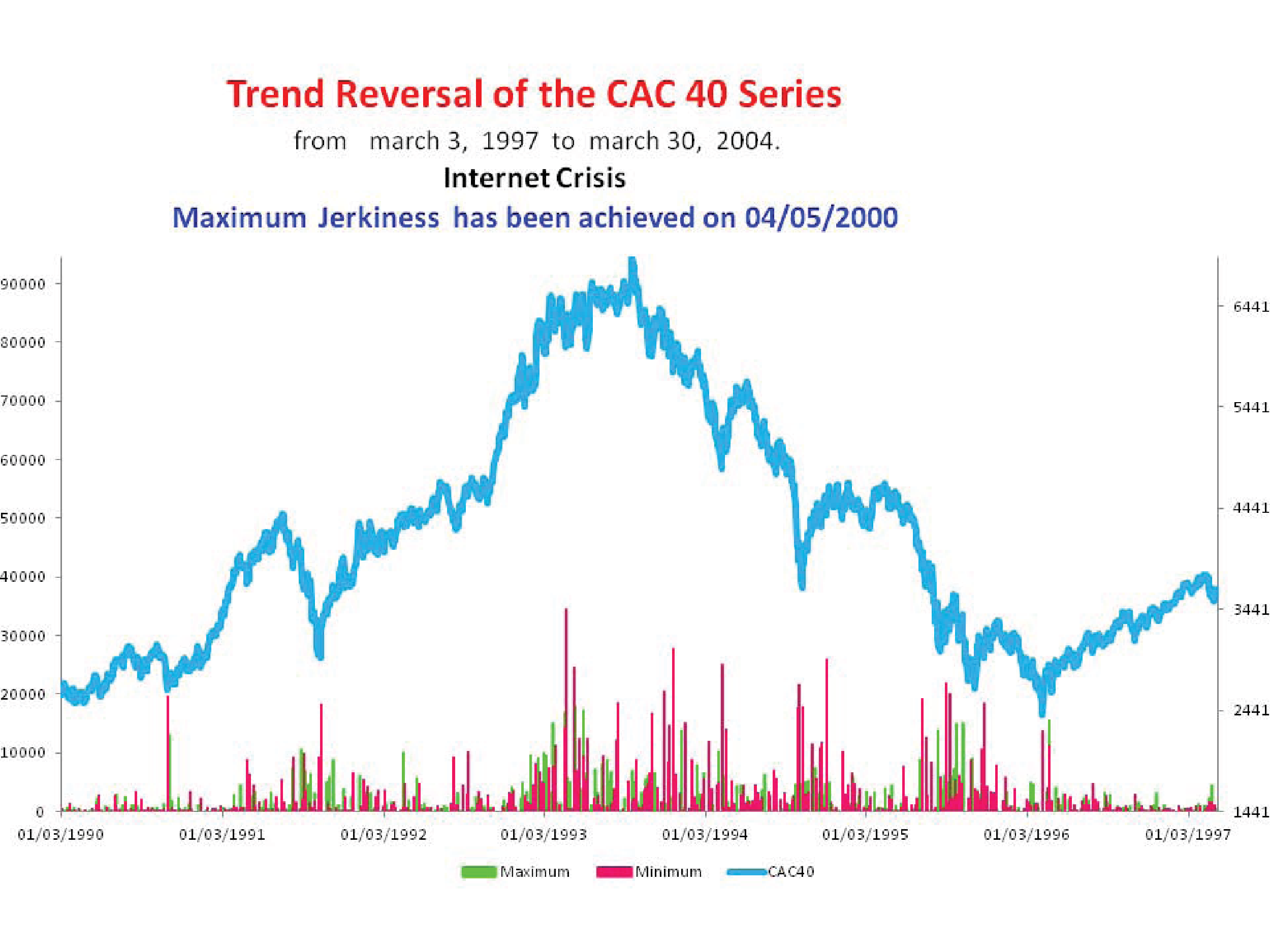}  & 
\includegraphics[width=\linewidth]{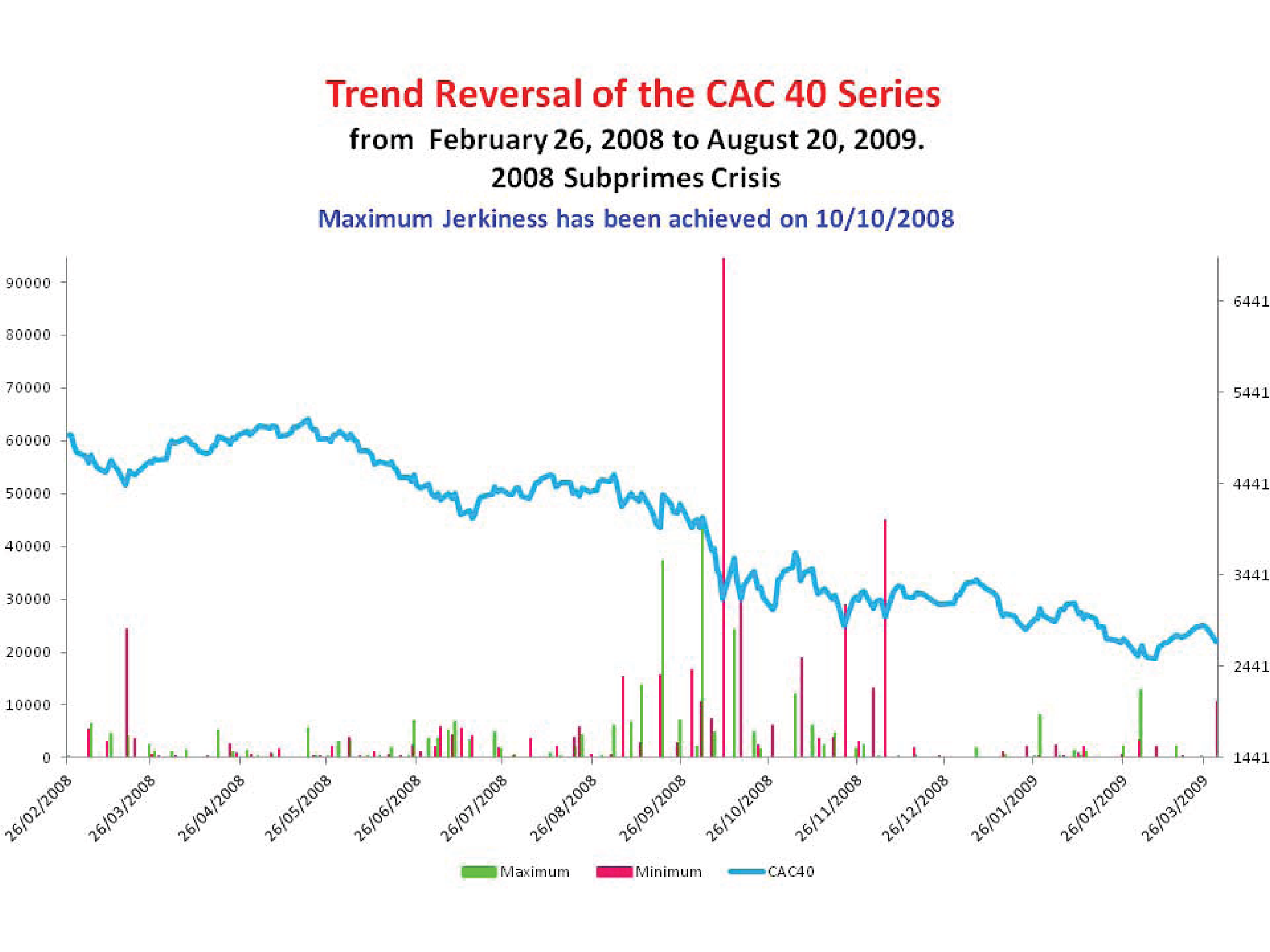}\\ \end{tabular} 
\end{center}\normalsize

The next figure displays the velocities of the jerkiness between 
two consecutive trend reversal dates, a ratio involving the 
variation of the jerkiness and the duration of the congruence 
period (bull and bear): 
 
 \vspace{ -1 mm}

\begin{center} 
\includegraphics[width=.7\linewidth]{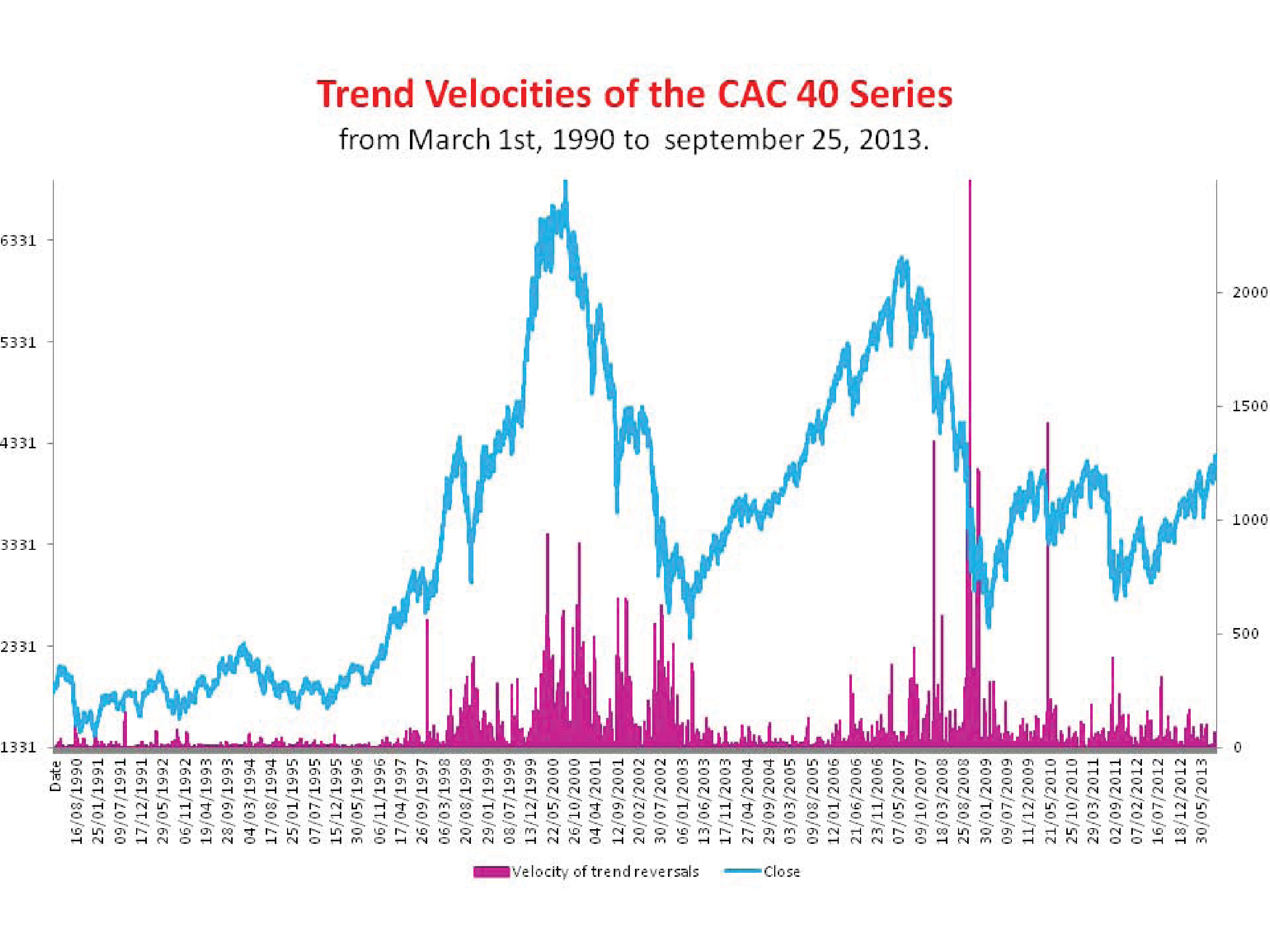} 
\vfill \mbox{} 
\end{center}
The following one displays the classification of trend speeds and 
absolute value of the accelerations  by decreasing jerkiness: 
\begin{center}\begin{tabular}{m{0.474\textwidth}m{0.474\textwidth}} 
\includegraphics[width=\linewidth]{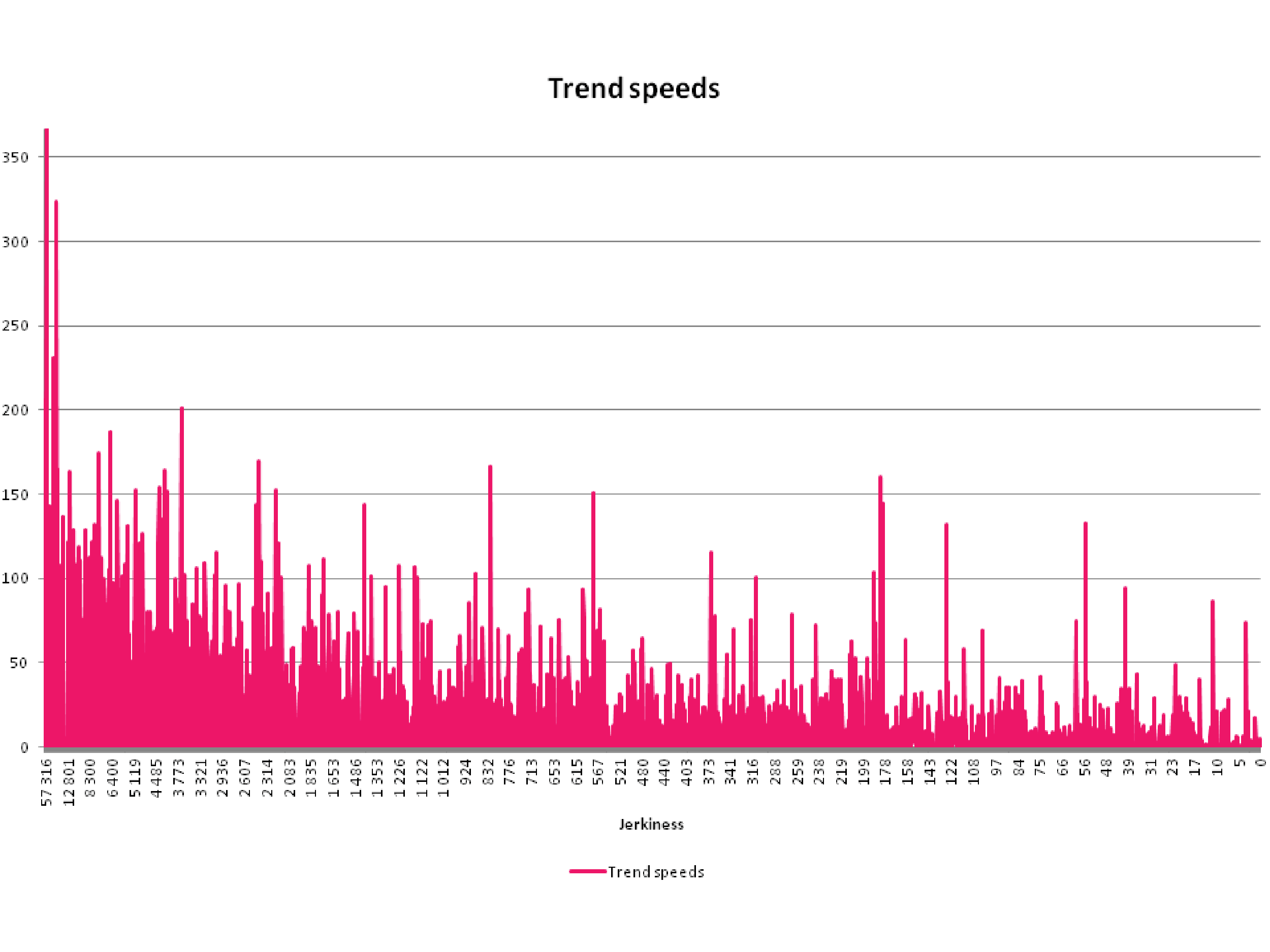}  & 
\includegraphics[width=\linewidth]{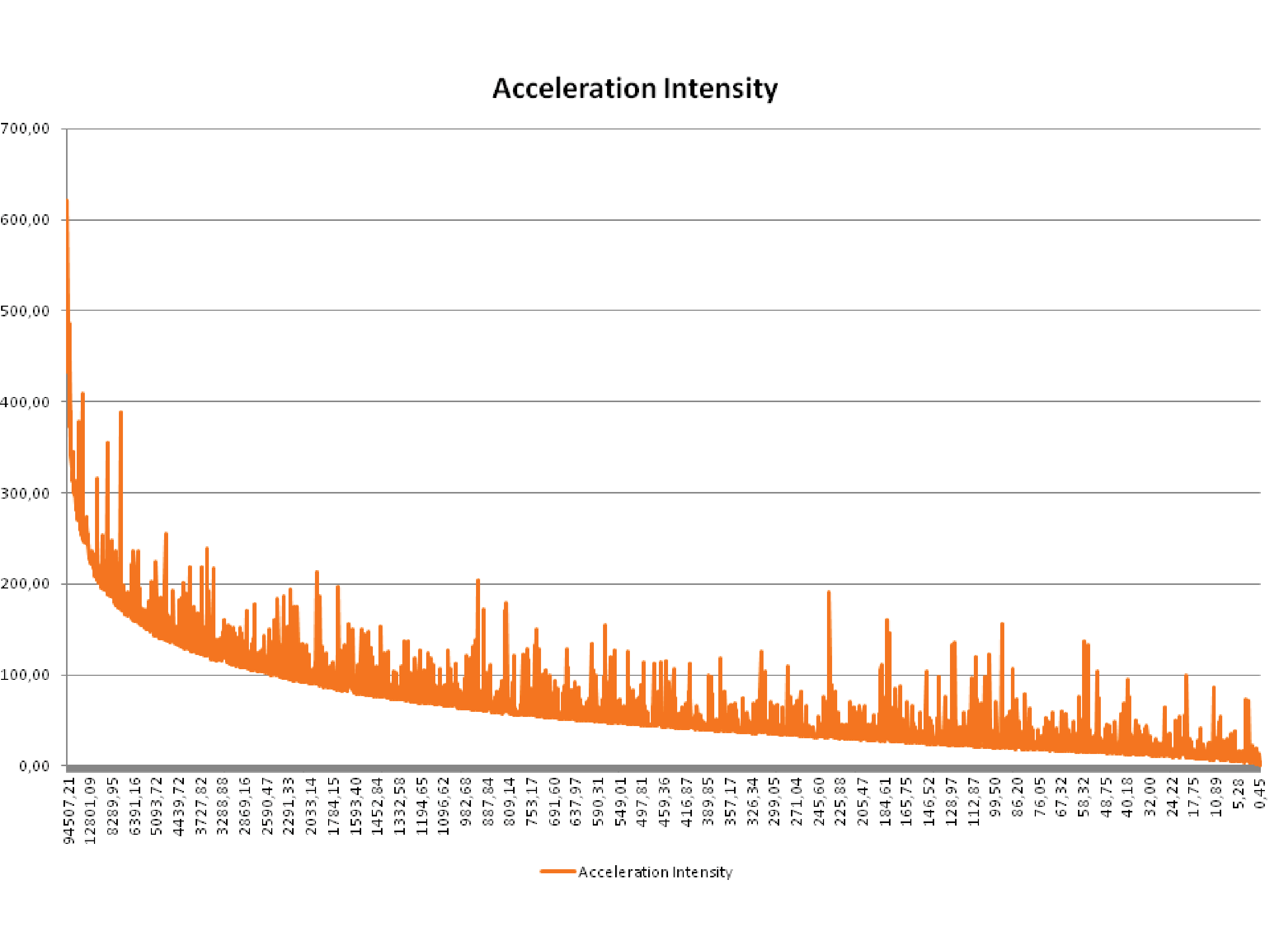}\\ 
\end{tabular}  \end{center}\normalsize

The analysis of this series  shows that often  the jerkiness at 
minima (bear periods) is higher than the ones at maxima (bull 
periods). For the CAC 40, the proportion of ``bear jerkiness'' 
(57\%) over ``bull jerkiness'' (43\%).

The next table provides the first dates by decreasing jerkiness.  
The most violent are those of the subprime crisis (in bold), then 
the ones of the year 2006 and, next, the dates of the Internet 
crisis (in italics).

 \begin{center} 
\begin{tabular}{|c|c|c|c|c|c|}  
\hline
\emph{\textbf{Date}}    &   \textbf{Jerkiness}   &   \emph{\textbf{Date}}    &   \textbf{Jerkiness}   &   \emph{\textbf{Date }}   &   \textbf{Jerkiness}   \\
\hline
10/10/\textbf{\textbf{2008}}  &   94507,21    &   03/01/\emph{2001}  &   15153,31    &   17/02/\emph{2000}  &   10025,57    \\
23/01/\textbf{2008}  &   57315,90    &   11/09/\emph{2002}  &   15111,43    &   28/10/\emph{2002}  &   9962,69 \\
07/05/2010  &   53585,50    &   10/03/\emph{2000}  &   15055,45    &   01/09/1998  &   9917,22 \\
05/12/\textbf{2008}  &   44927,23    &   10/08/2011  &   15011,24    &   15/02/\textbf{2008}  &   9905,51 \\
03/10/\textbf{2008}  &   43319,41    &   27/08/\emph{2002}  &   14958,41    &   19/04/1999  &   9887,67 \\
19/09/\textbf{2008}  &   37200,13    &   22/11/\emph{2000}  &   14768,91    &   26/10/\emph{2001}  &   9556,17 \\
05/04/\emph{2000}  &   34609,80    &   03/04/\emph{2000}  &   14280,35    &   29/06/\emph{2000}  &   9470,44 \\
21/01/\textbf{2008}  &   34130,42    &   03/04/\emph{2001}  &   14003,47    &   25/02/\emph{2000}  &   9438,07 \\
16/10/\textbf{2008}  &   29794,42    &   18/07/\emph{2002}  &   13813,67    &   27/03/\emph{2001}  &   9436,84 \\
21/11/\textbf{2008}  &   28840,69    &   19/12/\emph{2000}  &   13743,01    &   15/05/\emph{2000}  &   9411,84 \\
04/12/\emph{\emph{2000}} &   27861,03    &   12/03/2003  &   13707,93    &   04/10/2011  &   9409,14 \\
12/11/\emph{2001}  &   26039,07    &   12/09/\textbf{2008}  &   13682,85    &   17/01/\emph{2000}  &   9398,39 \\
22/03/\emph{2001}  &   25128,11    &   01/12/\textbf{2008}  &   13207,66    &   11/08/1998  &   9320,83 \\
27/04/\emph{2000}  &   24577,70    &   29/10/1997  &   13085,95    &   20/11/\textbf{2007}  &   9291,91 \\
17/03/\textbf{2008}  &   24416,22    &   04/03/2009  &   12845,84    &   05/10/1998  &   9277,96 \\
14/10/\textbf{2008}  &   24007,60    &   14/03/\textbf{2007}  &   12801,09    &   29/07/1999  &   9253,97 \\
05/08/\emph{2002}  &   22021,61    &   24/06/\emph{2002}  &   12658,98    &   04/12/\textbf{2007}  &   9200,48 \\
14/09/\emph{2001}  &   21658,15    &   02/08/2012  &   12628,14    &   04/02/\emph{2000}  &   9093,25 \\
10/08/\textbf{2007}  &   21252,50    &   24/05/\emph{2000}  &   12456,94    &   02/10/\emph{2002}  &   8959,94 \\
13/11/\emph{2000}  &   20662,32    &   10/05/\emph{2000}  &   12411,27    &   13/09/\emph{2000}  &   8897,37 \\
22/01/\textbf{2008}  &   20184,96    &   28/07/\emph{2000}  &   12145,83    &   10/05/2010  &   8877,39 \\
14/08/\emph{2002}  &   20052,16    &   23/02/\emph{2001}  &   11960,59    &   30/09/\emph{2002}  &   8845,61 \\
28/10/1997  &   19720,61    &   04/11/\textbf{2008}  &   11904,50    &   04/11/1998  &   8843,75 \\
14/06/\emph{2002}  &   19114,56    &   08/06/\textbf{\emph{\textbf{\emph{2006}}}}  &   11773,65    &   09/08/2011  &   8833,20 \\
06/11/\textbf{2008}  &   18900,51    &   30/10/\emph{2001}  &   11733,86    &   11/06/\emph{2002}  &   8832,22 \\
03/08/\emph{2000}  &   18621,37    &   15/10/\emph{2001}  &   11630,50    &   07/07/\emph{2000}  &   8797,60 \\
29/10/\emph{2002}  &   18550,19    &   24/03/2003  &   11294,44    &   16/01/\emph{2001}  &   8778,74 \\
08/10/1998  &   18307,12    &   15/03/\emph{2000}  &   11232,52    &   27/04/1998  &   8721,52 \\
02/05/\emph{2000}  &   18087,38    &   17/09/\textbf{2007}  &   10948,51    &   19/02/\textbf{2008}  &   8327,20 \\
21/09/\emph{2001}  &   17771,78    &   13/08/\textbf{2007}  &   10933,30    &   20/11/\emph{2000}  &   8299,90 \\
11/09/\emph{2001}  &   17660,69    &   25/10/\emph{2001}  &   10809,42    &   03/07/\emph{2002}  &   8289,95 \\
16/08/\textbf{2007}  &   17398,86    &   02/10/\textbf{2008}  &   10720,31    &   28/06/\emph{2000}  &   8258,67 \\
16/05/\emph{2000}  &   17228,62    &   23/10/\emph{2002}  &   10675,86    &   28/06/2010  &   8137,05 \\
04/04/\emph{2000}  &   16958,95    &   25/08/1998  &   10673,02    &   31/01/\emph{2000}  &   8093,58 \\
18/10/\emph{2000}  &   16761,07    &   30/03/2009  &   10672,64    &   21/11/\emph{2000}  &   8074,23 \\
29/09/\textbf{2008}  &   16502,34    &   24/01/\textbf{2008}  &   10352,96    &   28/01/2009  &   8049,26 \\
08/08/\textbf{2007}  &   16048,09    &   20/03/\emph{2001}  &   10294,67    &   26/02/\textbf{2007}  &   8038,76 \\
21/03/2003  &   15703,11    &   14/12/\emph{2001}  &   10253,40    &   31/01/\emph{2001}  &   8033,95 \\
18/09/\textbf{2008}  &   15506,17    &   31/07/\textbf{2007}  &   10134,80    &   26/11/\emph{2002}  &   7933,90 \\
22/05/\textbf{\emph{2006}}  &   15470,19    &   26/04/\emph{2000}  &   10093,65    &   08/08/2011  &   7821,87 \\
05/09/\textbf{2008}  &   15406,87    &   02/09/1999  &   10080,12    &   18/05/2010  &   7793,80 \\
\hline\end{tabular}
\end{center} 
\clearpage

\clearpage  

\section{Other Examples of Differential Connection Tensors} \label{s:Remarks}

 
\subsection{Prospective and Retrospective Derivatives of Set-Valued 
Maps} \label{s:PRDerivatives}
 
We summarize the concept of graphical derivatives.  
\begin{Definition} 
\symbol{91}\textbf{Retrospective and Prospective Graphical 
Derivatives}\symbol{93} \label{d:PRderiv}\index{} Consider a 
set-valued map $F: X  \leadsto Y$ from a finite dimensional 
vector space $X$ to another one $Y$. Let  $(x,y) \in \mbox{\rm 
Graph}(F)$ an element  of its graph. We denote in this study by   
\begin{enumerate}   

\item  \emph{retrospective derivative}  $\overleftarrow{D}F(x,y): 
X  \leadsto Y$  associating with any direction $u \in X$ the set 
of elements $v \in   Y$ satisfying

\begin{equation} \label{e:}   
\liminf_{h \mapsto 0+, u_{h} \mapsto u} d \left(v, 
\frac{y-F(x-hu_{h}) }{h} \right) \; = \; 0 
\end{equation}
\item   \emph{prospective derivative} $\overrightarrow{D}F(x,y):X 
 \leadsto Y $ associating with any direction $u \in X$ the set of elements 
$v \in  Y$ satisfying

\begin{equation} \label{e:}   
\liminf_{h \mapsto 0+, u_{h} \mapsto u} d \left(v, 
\frac{F(x+hu_{h})-y}{h} \right) \; = \; 0 
\end{equation}
\end{enumerate}

\end{Definition}

The retrospective and prospective difference quotients of $F$ at 
$(x,y) \in \mbox{\rm Graph}(F)$ are defined by $\displaystyle{ 
\overleftarrow{\nabla }_{h}F(x,y)(\overleftarrow{u}) := 
\frac{y-F(x-h\overleftarrow{u})}{h}}$ and
 $\displaystyle{\overrightarrow{\nabla 
}_{h}F(x,y)(\overrightarrow{u}) 
:=\frac{F(x+h\overrightarrow{u})-y}{h}}$. 

We can reformulate the definition of the (contingent) derivative 
by saying that it is the \emph{upper Painlevé-Kuratowski limit} 
of the difference quotients, 
\begin{equation} \label{e:}   
\forall \; \overleftarrow{u}, \; \; 
\overleftarrow{D}F(x,y)(\overleftarrow{u}) \; = \; \mbox{\rm 
Limsup}_{h \mapsto 0+, u_{h} \rightarrow \overleftarrow{u}} 
\overleftarrow{\nabla }_{h}F(x,y)(u_{h})
\end{equation}
i.e., the retrospective (resp. prospective) derivatives are the 
cluster points $\overleftarrow{v}$ of 
$\displaystyle{\overleftarrow{v}_{h} \in   \overrightarrow{\nabla 
}_{h}F(x,y)(u_{h}) }$ (resp. of i.e.,   the cluster points
 of $\displaystyle{\overrightarrow{v}_{h} \in   
\overrightarrow{\nabla }_{h}F(x,y)(u_{h})}$).  Whenever the 
set-valued map $F$ is Lipschitz, \emph{the retrospective and 
prospective difference quotients   are bounded}, and thus, 
relatively compact set since the dimension of the vector spaces 
is finite. In this case, the prospective and retrospective 
derivatives are not empty.

\vspace{ 5 mm} Taking the tensor product of  both the 
retrospective and prospective  derivatives allows us to define 
the differential connection tensor:
 

\begin{Definition} 
\symbol{91}\textbf{Differential Connection Tensor}\symbol{93} 
\label{d:DifferentialConnectionTensor} \index{} The 
\index{differential connection tensor} \emph{differential 
connection tensor} 
$\mathbf{a}_{F}(x,y)[(\overleftarrow{u},\overrightarrow{u}), 
(\overleftarrow{v},\overrightarrow{v})]$ of retrospective and 
prospective derivatives of $F$ at $(x,y)  \in \mbox{\rm 
Graph}(F)$ is defined by  

\begin{equation} \label{e:DifferentialConnectionTensor} \left\{ \begin{array}{l}  
\forall \; (\overleftarrow{u},\overrightarrow{u}), \; 
\overleftarrow{v} \in  
\overleftarrow{D}F(x,y)(\overleftarrow{u}),\; \overrightarrow{v} 
\in  \overrightarrow{D}F(x,y)(\overrightarrow{u}), \\ 
\displaystyle{ \mathbf{a}_{F}(x,y) 
[(\overleftarrow{u},\overrightarrow{u}),(\overleftarrow{v},\overrightarrow{v})] 
\; := \;   \overleftarrow{v} \otimes \overrightarrow{v} }
\end{array} \right. \end{equation}
\end{Definition}

\textbf{Remark} --- \hspace{ 2 mm} A normalized version of the 
differential connection tensor is defined by 

\begin{equation} \label{e:} \left\{ \begin{array}{l}  
\forall \; (\overleftarrow{u},\overrightarrow{u}), \; 
\overleftarrow{v} \in  
\overleftarrow{D}F(x,y)(\overleftarrow{u}),\; \overrightarrow{v} 
\in  \overrightarrow{D}F(x,y)(\overrightarrow{u}), \\ 
\displaystyle{ \mathbf{a}_{F}(x,y) 
[(\overleftarrow{u},\overrightarrow{u}),(\overleftarrow{v},\overrightarrow{v})] 
\; := \;  \frac{  \overleftarrow{v} \otimes \overrightarrow{v} }{ 
\|\overleftarrow{v}\|  \| \overrightarrow{v}\|}}
\end{array} \right. \end{equation}

The normalized version is not that useful whenever we are 
interested to the signs of the entries of the connection matrix. 
\hfill $\;\; \blacksquare$ \vspace{ 5 mm}

\textbf{Remark} --- \hspace{ 2 mm} One can associate with the  
prospective difference quotient  $\overrightarrow{\nabla 
}_{h}F(x,y)(\overrightarrow{u})$ and retrospective difference 
quotient $\overleftarrow{\nabla }_{h}F(x,y)(\overrightarrow{u})$ 
their difference quotient 
\begin{equation} \label{e:}  
\nabla ^{2}F(x,y)(\overleftarrow{u},\overrightarrow{u}) \; := \; 
\frac{\overrightarrow{\nabla }_{h}F(x,y)(\overrightarrow{u})- 
\overleftarrow{\nabla }_{h}F(x,y)(\overrightarrow{u}) }{h} \; = 
\; \frac{F(x+h\overrightarrow{u}) + F(x-h\overleftarrow{u})-2y 
}{h^{2}} 
\end{equation}
The Painlevé-Kuratowski upper limit of $\nabla 
^{2}F(x,y)(\overleftarrow{u},\overrightarrow{u})$ defines the 
retrospective-prospective second order graphical derivative of 
$F$ at $(x,y) \in \mbox{\rm Graph}(F)$ by:

\begin{equation} \label{e:}   
D^{2}F(x,y)(\overleftarrow{u},\overrightarrow{u}) \; := \; 
\mbox{\rm Limsup}_{h \mapsto 0+, \overleftarrow{u}_{h} 
\rightarrow \overleftarrow{u}, \overrightarrow{u}_{h} \rightarrow 
\overrightarrow{u} }  \nabla 
^{2}F(x,y)(\overleftarrow{u}_{h},\overrightarrow{u}_{h})
\end{equation}
The differential connection tensor replaces the difference 
between the retrospective and prospective derivatives by their 
tensor products.  We refer to Section~5.6, p. 315, of  
\emph{Set-valued analysis}, \cite[Aubin \& Frankowska]{af90sva}, 
for other approaches of higher order graphical derivatives of 
set-valued maps. \hfill $\;\; \blacksquare$ \vspace{ 5 mm} 

\textbf{Remark} --- \hspace{ 2 mm} In 1884,  \glossary{Peano 
(Giuseppe) [1858-1932]} \emph{Giuseppe Peano} proved in 
\glossary{Peano (Giuseppe) [1858-1932]} \emph{Giuseppe Peano} See 
\cite[\emph{Applicazioni geometriche del calcolo 
infinitesimale}]{Peano1887} that  continuous   derivatives are 
the limits 
\begin{displaymath}    \label{e:}   
\forall \; t \in \symbol{93}a,b\symbol{91}, \; \; \lim_{h 
\rightarrow 0}\frac{x(t+h)-x(t-h)}{2h} \; = \; 
 \frac{1}{2} \left(  \lim_{h \rightarrow 0+} 
\frac{x(t)-x(t-h)}{h}+ \lim_{h \rightarrow 
0}\frac{x(t+h)-x(t)}{h}\right)
\end{displaymath}
of  both the retrospective and prospective average velocities 
(difference quotients) at time $t$. We follow his suggestion  by 
taking the average of the  prospective difference quotient  
$\overrightarrow{\nabla }_{h}F(x,y)(\overrightarrow{u})$ and 
retrospective difference quotient $\overleftarrow{\nabla 
}_{h}F(x,y)(\overleftarrow{u})$ their difference quotient 
\begin{equation} \label{e:}  
\frac{\overrightarrow{\nabla }_{2h}F(x,y)(\overrightarrow{u})+ 
\overleftarrow{\nabla }_{h}F(x,y)(\overleftarrow{u}) }{2h}  
\end{equation} 
and taking their Painlevé-Kuratowski limits
\begin{equation} \label{e:}  
\mbox{Limsup}_{h \mapsto 0+, \overrightarrow{u}_{h} \rightarrow 
\overrightarrow{u}}  \overrightarrow{\nabla} 
_{h}F(x,y)(\overrightarrow{u}_{h})+ \mbox{Limsup}_{h \mapsto 0+, 
\overleftarrow{u}_{h} \rightarrow \overleftarrow{u}} 
\overleftarrow{\nabla }_{h}F(x,y)(\overleftarrow{u}_{h})    
\end{equation} 
in order to define \index{Peano graphical derivative} \emph{Peano 
graphical derivatives} of $F$ at $(x,y) \in \mbox{\rm Graph}(F)$ 
depending on \emph{pairs $(\overleftarrow{u},\overrightarrow{u})$ 
of directions}.   \hfill $\;\; \blacksquare$ \vspace{ 5 mm}

 
\subsection{Differential Connection Tensors of Numerical 
Functions}  \label{s:PREpider}

When $V: x \in X \mapsto V(x) \in \{-\infty \} \cup\mathbb{R}\cup 
\{+\infty \}$ is an extended numerical function on 
$\mathbb{R}^{}$, it can also be regarded as a set-valued map 
(again denoted by) $V: X  \leadsto \mathbb{R}^{}$ defined by 
\begin{equation} \label{e:} V(x) \; := \left\{ \begin{array}{ccc}
  \{V(x)\} &\mbox{\rm if}& V(x) \in \mathbb{R}^{} \; \;(\mbox{\rm 
 i.e.,} \; x \; \in \; \mbox{\rm Dom}(V)) \\
 \emptyset  &\mbox{\rm if not}& 
 \end{array} \right. \end{equation} 

A slight modification of  Theorem~6.1.6, p. 230 of 
\emph{Set-valued analysis}, \cite[Aubin \& Frankowska]{af90sva}, 
states that 

\begin{equation} \label{e:} \left\{ \begin{array}{l}   
\overrightarrow{D}V(x)(\overrightarrow{u})\;  = \; 
[\overrightarrow{D}_{\uparrow}V(x)(\overrightarrow{u}), 
\overrightarrow{D}_{\downarrow} V(x)(\overrightarrow{u})] \\ 
\mbox{} \\ \overleftarrow{D}V(x)(\overleftarrow{u})\;  = \; 
[\overleftarrow{D}_{\uparrow}V(x)(\overleftarrow{u}), 
\overleftarrow{D}_{\downarrow} V(x)(\overleftarrow{u})]
\end{array} \right. 
\end{equation} 
where

\begin{equation} \label{e:} \left\{ \begin{array}{l}  
\displaystyle{\overrightarrow{D}_{\uparrow}V(x)(\overrightarrow{u}) 
\; := \; \liminf_{h \rightarrow 0+} 
\frac{V(x+h\overrightarrow{u})-V(x)}{h} \;\mbox{\rm 
(epiderivative of $V$)}\;
}\\
\displaystyle{\overrightarrow{D}_{\downarrow} 
V(x)(\overrightarrow{u})\; := \; \limsup_{h \rightarrow 0+} 
\frac{V(x+h\overrightarrow{u})-V(x)}{h} \;\mbox{\rm 
(hypoderivative of $V$)}\;
}\\
 \displaystyle{\overleftarrow{D}_{\uparrow}V(x)(\overleftarrow{u}) 
\; := \; \liminf_{h \rightarrow 0+} \frac{V(x 
)-V(x-h\overleftarrow{u})}{h} \; = \; 
- \overrightarrow{D}_{\downarrow} V(x)(-\overleftarrow{u})}\\
\displaystyle{\overleftarrow{D}_{\downarrow} 
V(x)(\overleftarrow{u})  \; := \; \limsup_{h \rightarrow 0+} 
\frac{V(x )-V(x-h\overleftarrow{u})}{h} \; = \; - 
\overrightarrow{D}_{\uparrow} V(x)(-\overleftarrow{u})
}\\
\end{array} \right. 
\end{equation}

Definition~\ref{d:DifferentialConnectionTensor}, 
p.\pageref{d:DifferentialConnectionTensor} implies that  
\begin{equation} \label{e:DifferentialConnectionTensor} \left\{ \begin{array}{l}  
\forall \; (\overleftarrow{u},\overrightarrow{u}), \; 
\overleftarrow{v} \in  
\overleftarrow{D}V(x)(\overleftarrow{u}),\; \overrightarrow{v} 
\in  \overrightarrow{D}V(x)(\overrightarrow{u}), \\ 
\displaystyle{ \mathbf{a}_{V}(x,y) 
[(\overleftarrow{u},\overrightarrow{u}),(\overleftarrow{v},\overrightarrow{v})] 
\; := \;   \overleftarrow{v}   \overrightarrow{v} }
\end{array} \right. \end{equation}
since tensor products of real numbers boil down to their 
multiplication.

Therefore, for any pair $(\overleftarrow{u},\overrightarrow{u})$, 
the subset of differential connection tensors of retrospective 
and prospective directions is equal to

\begin{equation} \label{e:} \left\{ \begin{array}{l}
\overleftarrow{D}V(x)(\overleftarrow{u}) \otimes 
\overrightarrow{D}V(x)(\overrightarrow{u}) \; := \\ 
\displaystyle{\left\{ \overleftarrow{v} \overrightarrow{v} 
\right\} _{(\overleftarrow{v},\overrightarrow{v}) \in 
[\overleftarrow{D}_{\uparrow}V(x)(\overleftarrow{u}), 
\overleftarrow{D}_{\downarrow} V(x)(\overleftarrow{u})] \times 
[\overrightarrow{D}_{\uparrow}V(x)(\overrightarrow{u}), 
\overrightarrow{D}_{\downarrow} V(x)(\overrightarrow{u})] }}
\end{array} \right. 
\end{equation}

\begin{Definition} 
\symbol{91}\textbf{Reversal Direction Pair}\symbol{93} 
\label{d:RevesDirPair}\index{} A pair 
$(\overleftarrow{u},\overrightarrow{u}) $ of directions 
$\overleftarrow{u} \in X$ and $\overrightarrow{u} \in X $ is a 
\index{reversal direction pair} \emph{reversal direction pair} of 
$V$ at $x \in \mbox{\rm Dom}(V)$ if 

\begin{equation} \label{e:}   
\overleftarrow{D}_{\uparrow}(x)(\overleftarrow{u}) 
\overrightarrow{D}(x)(\overrightarrow{u}) \; = \; 
\overleftarrow{D}_{\downarrow} (x)(-\overrightarrow{u} )
\overrightarrow{D}_{\downarrow} V(x)(-\overleftarrow{u}) \; < \; 0
\end{equation}
A direction $u \in X$ is a reversal direction of $V$ at $x$ if 
the diagonal pair $(u,u)$ is reversal  direction pair. \\ This 
means that a positive (resp. negative) retrospective 
epiderivative of $V$ at $x$ in the direction $\overleftarrow{u}$ 
is followed  by a negative (resp. positive) prospective 
epiderivative in the direction $\overrightarrow{u}$, or, 
respectively,that a positive (resp. negative) retrospective 
hypoderivative in the direction $- \overrightarrow{u}$ is 
followed by a negative (resp. positive) prospective 
hypoderivative in the direction $-\overleftarrow{u}$. 
\end{Definition}

Recall that if $V$ achieves a local minimum at $x$, the Fermat 
rule states  that

\begin{equation} \label{e:}   
\forall \; \overrightarrow{u} \in X, \; \;  
\overrightarrow{D}_{\uparrow}V(x)(\overrightarrow{u})   \; \geq 
\; 0 \;\mbox{\rm and}\; \forall \; \overleftarrow{u} \in X, \; 
\;  \overleftarrow{D}_{\downarrow} V(x)(\overleftarrow{u})  \; 
\leq \; 0
\end{equation}
and if it achieves a local maximum at $x$, that
\begin{equation} \label{e:}   
\forall \; \overrightarrow{u} \in X, \; \;  
\overrightarrow{D}_{\downarrow} V(x)(\overrightarrow{u})   \; 
\leq \; 0 \;\mbox{\rm and}\; \forall \; \overleftarrow{u} \in X, 
\; \;  \overleftarrow{D}_{\uparrow} V(x)(\overleftarrow{u})  \; 
\geq \; 0
\end{equation}
 
These conditions are not sufficient for characterizing local 
extrema: convexity or many second order conditions provide 
sufficient conditions (see \emph{Set-valued 
analysis},\cite[Aubin  \& Frankowska]{af90sva}, \emph{Variational 
Analysis}, \cite[Rockafellar \& Wets]{rw91nsa} and an important 
literature on set-valued and variational analysis).

Recall that the prospective epidifferential (or prospective 
epidifferential subdifferential) $\overrightarrow{\partial} 
_{\uparrow}V(x)$ of a function $V$ at $x$ is the set of elements 
$\overrightarrow{p}_{\uparrow} \in X^{\star}$ such that for any 
$v \in X$, $\left\langle  \overrightarrow{p}_{\uparrow},v 
\right\rangle  \; \leq \; \overrightarrow{D}_{\uparrow}V(x)(v)$. 
In the same way, we define the retrospective epidifferential (or 
retrospective epidifferential subdifferential)  
$\overleftarrow{\partial} _{\uparrow}V(x)$ of a function $V$ at 
$x$ as the set of elements $\overleftarrow{p}_{\uparrow} \in 
X^{\star}$ such that for any $v \in X$, $\left\langle  
\overleftarrow{p}_{\uparrow},v \right\rangle  \; \leq \; 
\overleftarrow{D}_{\uparrow}V(x)(v)$.  It is equal to prospective 
hypodifferential (or prospective superdifferential) 
$\overrightarrow{\partial}_{\downarrow} V(x) $,  the set of 
elements $\overrightarrow{p}_{\downarrow}  \in X^{\star}$ such 
that for any $v \in X$, $\left\langle 
\overrightarrow{p}_{\downarrow},v \right\rangle  \; \geq \; 
\overrightarrow{D}_{\downarrow} V(x)(v)$. 

For instance, the trendometer detects the local extrema of 
numerical functions, such as the function $t \mapsto 
1-cos(2t)cos(3t)$: 
 
\begin{center} \includegraphics[width=.8\linewidth]{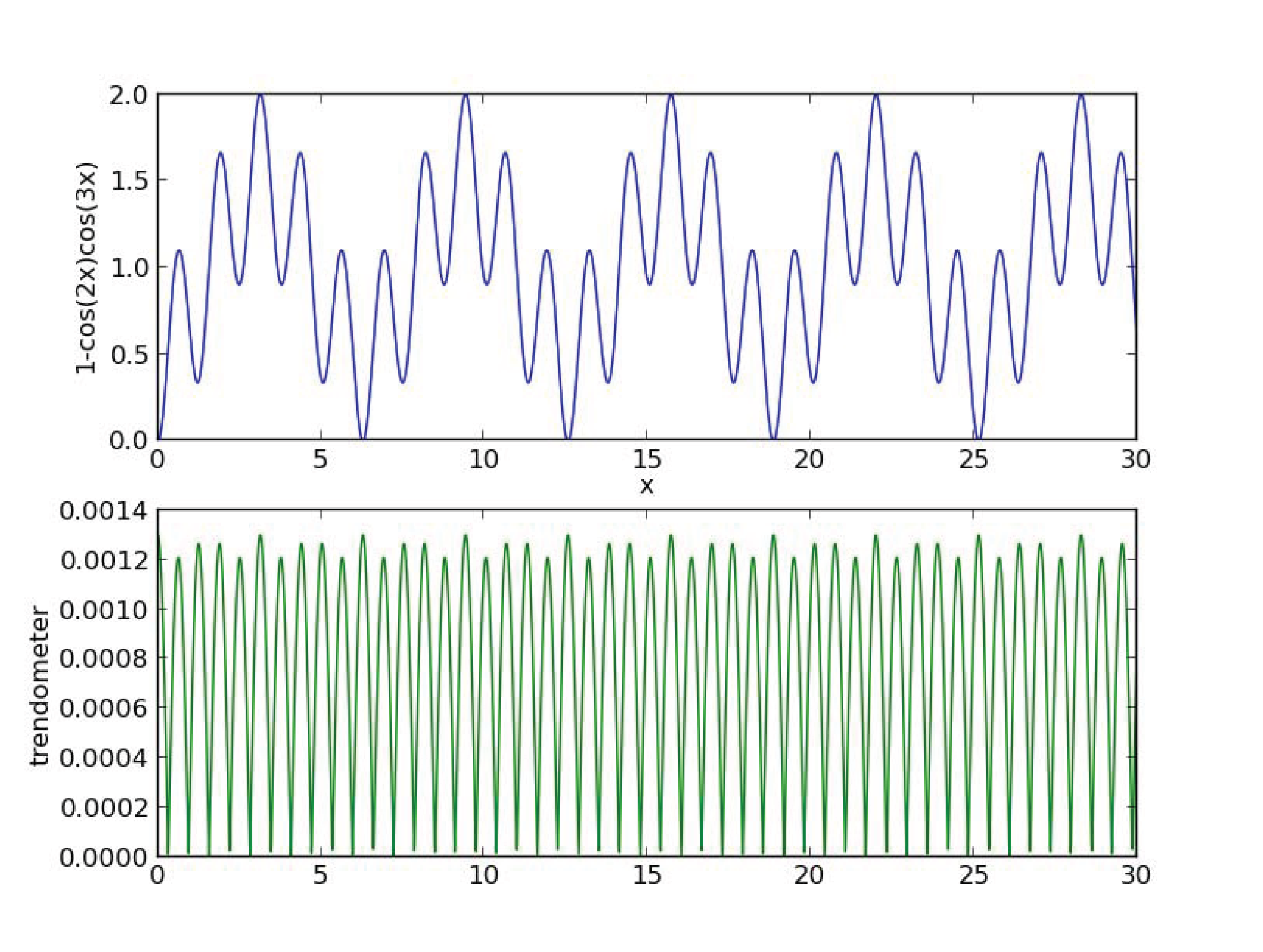}  
 \end{center}
 
\subsection{Tangential Connection Tensors} \label{s:PRCones}

The tangent spaces to differentiable manifolds being vector 
spaces, directions arriving at a point (we may call them 
\index{retrospective tangent} \emph{retrospective})  and 
directions starting from this point  \index{prospective tangent} 
(\emph{prospective}) belong to the same vector space. This is no 
longer the case when the subset is any (closed) subset $K \subset 
X$ of a finite dimensional vector space $X$. However, we may 
replace  vector spaces by cones.

We are indebted to the historical studies \cite[Dolecki \& 
Greco]{DoleckiGreco} (in which the authors quote 
\glossary{Fréchet (Maurice) [1878-1973]} \emph{Maurice Fréchet} 
stating that ``Cette théorie des ``contingents et paratingents" 
dont l'utilité a été signalée d'abord par M. Beppo Levi, puis par 
M. Severi, mais dont on doit à M. Bouligand et ses élèves d'en 
avoir entrepris l'étude systématique.'') and \cite[Greco, 
Mazzucchi \&  Pagani]{GrecoMazzucchiPagani}.  \glossary{Severi 
(Francesco) [1879-1961]} \emph{Francesco Severi} and 
\glossary{Bouligand (Georges) [1889-1979]} \emph{Georges 
Bouligand}, a whole menagerie of tangent cones, the definitions 
of which depend upon the limiting process, have been proposed 
(among many monographs, see \cite[\emph{Set-valued 
analysis}]{af90sva} and  \emph{Variational Analysis}, 
\cite[Rockafellar \& Wets]{rw91nsa} for instance). At some 
points, the tangent cones are not vector spaces, and the opposite 
of some tangent directions may no longer be tangent. 
 
We  suggest to regard the (contingent) tangent cone\footnote{See 
\cite[\emph{Set-valued analysis}]{af90sva}. The  (adjacent)
Peano-Severi-Bouligand tangent cone is defined by the 
Painlevé-Kuratowski lower limits instead of upper limits
\begin{equation} \label{e:}   
 \mbox{Liminf}_{h \mapsto 0+} 
\frac{K-x}{h} \; := \; \left\{ \overrightarrow{v} \; \in \; X \; 
\mbox{ such that} \; \lim_{h \mapsto 0+} 
\frac{d_{K}(x+h\overrightarrow{v})}{h} \; = \; 0\right\}
\end{equation}
The smaller  \emph{adjacent}   tangent cone is used whenever more 
regularity is required. An element $x \in K$ is said to be 
\index{regular element in a subset} \emph{regular} in $K$ at $x$ 
if both contingent and adjacent tangent cones coincide, i.e., 
when $T_{K}(x)$ is the Painlevé-Kuratowski limit of 
$\displaystyle{\frac{K-x}{h}}$.}  as the  \index{prospective 
tangent cone} \emph{prospective tangent cone} to $K$ at $x \in K$ 
defined by the Painlevé-Kuratowski upper limits 

\begin{equation} \label{e:}   
\overrightarrow{T}_{K}(x) \; := \; \mbox{Limsup}_{h \mapsto 0+} 
\frac{K-x}{h} \; := \; \left\{ \overrightarrow{v} \; \in \; X \; 
\mbox{ such that} \; \liminf_{h \mapsto 0+} 
\frac{d_{K}(x+h\overrightarrow{v})}{h} \; = \; 0\right\}
\end{equation}
with which we associate \index{retrospective tangent cone} 
(adjacent) \emph{retrospective tangent cone}\footnote{Backward 
evolutions and negative tangents have been introduced in 
\cite[Frankowska]{fh91cdc,HJB92} for characterizing lower 
semicontinuous (viscosity) solutions to Hamilton-Jacobi-Bellman 
equations.} 

\begin{equation} \label{e:}   
\overleftarrow{T}_{K}(x) \; := \; \mbox{Limsup}_{h \mapsto 0+} 
\frac{x-K}{h} \; :=   \; := \;  \left\{ \overleftarrow{v} \; \in 
\; X \; \mbox{ such that} \; \liminf_{h \mapsto 0+} 
\frac{d_{K}(x-h\overleftarrow{v})}{h} \; = \;  \right\}
\end{equation}  
satisfying $\overleftarrow{T}_{K}(x) \; := \; 
-\overleftarrow{T}_{K}(x)$.     It is natural to consider their 
tensor product $(x-h\overleftarrow{v}) \otimes  (x+ h 
\overrightarrow{v})$.   
 The  signs of its  entries   detect the ``blunt'' and``sharp''   elements of the 
boundary  in the same directions  \index{trend congruence} 
(\emph{trend congruence}) or in opposite directions \index{trend 
reversal} (\emph{trend reversal}).

 \clearpage

\tableofcontents

\begin{thebibliography}{a}
 
  



\bibitem{a83slhea} Aubin J.-P.    (1983)
Slow and heavy trajectories of controlled problems: smooth 
viability domains, In \emph{Multifunctions and Integrands}, 
105-116 Lecture Notes in Mathematics,\#1091, Ed. Salinetti G., 
Springer-Verlag 

\bibitem{avt} Aubin J.-P.    (1991)
 \emph{Viability Theory}, Birkh\"auser

 

\bibitem{a92ia} Aubin J.-P.    (1996)
\emph{Neural Networks and Qualitative Physics: a Viability 
Approach}, Cambridge University Press

\bibitem{aconcom98} Aubin J.-P.    (1998)
Connectionist Complexity and its Evolution, in \emph{Equations 
aux d\'eriv\'ees partielles, Articles d\'edi\'es \`a J.-L. 
Lions}, 50-79, Elsevier 
 
\bibitem{a01ctfcr} Aubin J.-P.    (2003)
Regulation of the Evolution of the Architecture of a Network by 
Connectionist Tensors Operating on Coalitions of Actors, \emph{J. 
Evolutionary Economics}, 13,95-124

 

\bibitem{AUB-Leit09}  Aubin J.-P. (2010)
Macroscopic Traffic Models: Shifting from Densities to 
``Celerities'', \emph{Applied Mathematics and Computation}, 217, 
963-971, \url{http://dx.doi.org/10.1016/j.amc.2010.02.032}


\bibitem{aub-12-Durance} Aubin J.-P. (2013)
Chaperoning State Evolutions by Variable Durations, \emph{SIAM 
Journal of Control and Optimization}, DOI. 10.1137/120879853 



\bibitem{TransReg-2013} Aubin J.-P.  (submitted)
Transports Regulators of  Networks with Junctions Detected by 
Durations Functions,


\bibitem{absp} Aubin J.-P., Bayen A. \& Saint-Pierre P. (2011)
 \emph{Viability Theory.  New Directions}, Springer

\bibitem{ab99nng}  Aubin J.-P.  \& Burnod Y.    (1998)
Hebbian Learning in Neural Networks with Gates,  \emph{Cahiers du 
Centre de Recherche Viabilit\'e, Jeux, Contrôle} \# 981

\bibitem{ACD-ALIM} Aubin J.-P., Chen Lx \& Dordan O. (2014) 
\emph{Tychastic Measure of Viability Risk. A 
Viabilist Portfolio Performance and Insurance Approach}

\bibitem{af90sva} Aubin J.-P. \& Frankowska H. (1990) 
 \emph{Set-valued analysis}, Birkhäuser
 
\bibitem{af85heav} Aubin J.-P. \& Frankowska H.   (1985)
 Heavy viable trajectories of controlled systems, 
Proceedings of  \emph{Dynamics of Macrosystems},  IIASA, 
September 1984,, Ed. Aubin J.-P., Saari D.\& Sigmund K., 
Springer-Verlag,148-167

\bibitem{ah01hyb} Aubin J.-P. \& Haddad G.   (2001)
Path-dependent impulse and hybrid control systems, in 
\emph{Hybrid Systems: Computation and Control}, 119-132, Di 
Benedetto \& Sangiovanni-Vincentelli Eds, Proceedings of the HSCC 
2001 Conference, LNCS 2034, Springer-Verlag

\bibitem{aubhadclio}   Aubin J.-P. \& Haddad G.    (2002) History
(Path) Dependent Optimal Control and Portfolio Valuation and 
Management, \emph{J. Positivity}, 6, 331-358
\bibitem{Buquoy}   Buquoy G. 
(1815) \emph{Exposition d'un nouveau principe général de 
dynamique, dont le principe des vitesses virtuelles n'est qu'un 
cas particulier}, V. Courtier

 
 
\bibitem{DoleckiGreco} 
Dolecki S. \& Greco G. H. (2007) Towards Historical Roots of 
Necessary Conditions of Optimality: Regula of Peano, 
\emph{Control and Cybernetics}, 36, 491-518    

\bibitem{do93liv} Dordan O. (1995) 
\textbf{\emph{Analyse qualitative}}, Masson



\bibitem{fh91cdc}
 Frankowska H. (1991) Lower semicontinuous solutions to 
Hamilton-Jacobi-Bellman equations, in \emph{Proceedings of the 
30th IEEE Conference on Decision and Control}, Brighton, UK, 


\bibitem{HJB92}
 Frankowska H. (1993) Lower semicontinuous solutions of 
Hamilton-Jacobi-Bellman equations, \emph{SIAM J. Control Optim.}, 
31,   257-272.


 
\bibitem{Galperin1} Galperin E. A. (2009) 
Information transmittal, time uncertainty, and special 
relativity, \emph{Computers and Mathematics with Applications}, 
  57,  1554-1573, doi: 
10.1016/j.camwa.2008.09.048

\bibitem{Galperin2} Galperin E. A. (2011)
Left time derivatives in mathematics, mechanics and control of 
motion, \emph{Computers and Mathematics with Application},   62
 4742–4757

\bibitem{Galperin3}
Galperin E. A. (submitted)  Information  Transmittal, Causality, 
Relativity  and Optimality,
  
\bibitem{Galperin4}
Galperin E. A. (submitted) Time And Relativity  In Dynamical 
Systems,
  
\bibitem{GrecoMazzucchiPagani} 
 Greco G. H.,   Mazzucchi S. \&  Pagani E. M. (2010) 
Peano on derivative of measures: strict derivative of 
distributive set functions, \emph{Rend. Lincei Mat. Appl.}, 21, 
305-339 DOI 10.4171/RLM/575  


\bibitem{hg81} Haddad G. (1981)
  Monotone trajectories of differential inclusions with 
memory, \emph{Isr. J. Math.}, 39, 83-100

\bibitem{hg81b} Haddad G. (1981)
Monotone viable trajectories for functional differential 
inclusions, \emph{J. Diff. Eq.}, 42, 1-24

\bibitem{hg81c} Haddad G. (1981)
Topological properties of the set of solutions for functional
  differential inclusions, \emph{Nonlinear Anal. Theory,
Meth. Appl.}, 5, 1349-1366


\bibitem{Hale} Hale J. K. (1993) 
\emph{Introduction to Functional Differential Equations}, Springer
 



\bibitem{Hebb} Hebb D. (1949) \emph{The Organization of Behavior}, Wiley
 
\bibitem{Mestschersky} Mestschersky I.V.  (1897) 
Dynamics of point with variable mass. In I.V. \emph{Mestschersky, 
Works on Mechanics of Bodies with Variable Mass},  Edition, 
Gostechizdat, Moscow, 1952,  37-188  

\bibitem{Levi-Civita} Levi-Civita (1928) Sul moto di un corpo di massa variablie,   
\emph{Rendiconti dei Lincei}, 329-333

\bibitem{Peano1887} Peano G. (1887) 
\emph{Applicazioni geometriche del calcolo infinitesimale}, 
Fratelli Bocca Editori, 
\url{http://historical.library.cornell.edu/cgi-bin/cul.math/docviewer?did=00610002&seq=1}

\bibitem{DanchinAll11} Porcar M., Danchin A. de Lorenzo V., dos Santos 
V., Krasnogor N., Rasmussen S.  \& Moya A. 
 (2011), The ten grand challenges of synthetic 
life, \emph{Syst Synth Biol}, 5,  1–9,  doi:  
10.1007/s11693-011-9084-5 


\bibitem{rw91nsa} Rockafellar R.T. \& Wets R. (1997) 
\emph{Variational Analysis}, Springer-Verlag


\bibitem{tp90vtb} Tallos P. (1991) Viability problems for nonautonomous 
differential inclusions, \emph{SIAM J. on Control and 
Optimization}, 29 
 253-263  
 
\bibitem{vino} Vinogradova G. (2012) 
Correction of Dynamical Network's Viability by Decentralization 
by Price, \emph{Complex Systems}, 21, 37-55 
\end{thebibliography}
\end{document}